\magnification 1200
\hoffset 0.5 true cm
\hsize 15 true cm
\overfullrule=0mm
\input amssym.def
\font \tit=cmbx6 scaled \magstep4
\font\bu=cmss8
\def\pn{\par\noindent}

\def\ZZ{\Bbb{Z}}
\def\QQ{\Bbb{Q}}

\def\CC{\Bbb{C}}
\def\GG{\Bbb{G}}

\def\du{\hbox{\bu v}}
\def\G{{\cal G}_{\rm mot}}

\def\L{{\rm Lie}\,}
\def\Hom{\rm Hom}
\def\End{\rm End}

\def\W{\rm W}
\def\Gr{\rm Gr}
\def\gr{\rm gr}
\def\spec{{\rm Spec}\,}
\def\sp{{\rm Sp}\,}

\def\ok{\overline k}
\def\gal{{\rm Gal}(\ok/k)}
\def\galm{{\cal GAL}(\ok/k)}

\null\vskip+2 true cm

\pn {\centerline {\tit Motivic Galois theory for motives of niveau $\leq 1$ }}

\vskip 1 true cm

\pn 
{\centerline {Cristiana Bertolin}}

\vskip 3 true cm

\pn {\bf Abstract.} Let $\cal T$ be a Tannakian category over a field $k$ of characteristic 0 and $\pi({\cal T})$ its fundamental group. In this paper we prove that there is a bijection between the $\otimes$-equivalence classes of Tannakian subcategories of ${\cal T}$ and the normal affine group sub-${\cal T}$-schemes of $\pi ({\cal T})$. 
\pn We apply this result to the Tannakian category ${\cal T}_1(k)$ generated by motives of niveau $\leq 1$ defined over $k$, whose fundamental group is called the motivic Galois group $\G({\cal T}_1(k))$ of motives of niveau $\leq 1$.
 We find four short exact sequences of affine group sub-${\cal T}_1(k)$-schemes of $\G({\cal T}_1(k))$, correlated one to each other by inclusions and projections.
Moreover, given a 1-motive $M$, we compute explicitly the biggest Tannakian subcategory of the one generated by $M$, whose fundamental group is commutative.
\vskip 0.5 true cm

\pn {\bf Mathematics Subject Classifications (2000).} 14L15, 18A22
\vskip 0.5 true cm

\pn {\bf Key words.} Fundamental group of Tannakian categories, Artin motives, 1-motives, motivic Galois groups.

\vskip 1 true cm

\pn {\bf Introduction.}
\vskip 1 true cm

\par Let $k$ be a field of characteristic 0 and $\ok$ its algebraic closure. 
Let $\cal T$ be a Tannakian category over $k$. 
The tensor product of $\cal T$ allows us to define the notion of Hopf algebras in
the category 
${\rm Ind}{\cal T}$ of Ind-objects of $\cal T$. The category of affine group 
$\cal T$-schemes is the opposite of the category of Hopf algebras in ${\rm Ind}{\cal T}$.
\pn The fundamental group $\pi ({\cal T})$ of $\cal T$
 is the affine group $\cal T$-scheme ${\rm Sp}(\Lambda),$ whose Hopf algebra
$ \Lambda$ satisfies the following universal property: for each object $X$ of ${\cal T}$, there exists a morphism
$ X \longrightarrow \Lambda \otimes X$ functorial in $X$. Those morphisms 
$\{ X \longrightarrow \Lambda \otimes X \}_{X \in {\cal T}}$ define {\it an action of the fundamental group $\pi ({\cal T})$ on each object of ${\cal T}$.}
\pn For each Tannakian subcategory ${\cal T}'$ of ${\cal T}$, let $H_{\cal T}({\cal T}')$ be the kernel of the faithfully flat morphism of affine group ${\cal T}$-schemes 
$I:\pi ({\cal T}) \longrightarrow i\pi ({\cal T}')$ corresponding to the inclusion functor $i:{\cal T}' \longrightarrow {\cal T} $. In particular we have the short exact sequence of affine group $\pi ({\cal T})$-schemes

$$ 0 \longrightarrow H_{\cal T}({\cal T}') \longrightarrow \pi ({\cal T})
\longrightarrow \pi ({\cal T}')  \longrightarrow 0. \leqno(0.1)$$

\pn In [6] 6.6, Deligne proves that the Tannakian category ${\cal T}'$ is equivalent, as tensor category, to the subcategory of ${\cal T}$ generated by those object on which the action of 
$ \pi ({\cal T})$ induces a trivial action of $H_{\cal T}({\cal T}')$. 
In particular, this implies that the fundamental group of $\pi ({\cal T}')$ 
of ${\cal T}'$ is isomorphic to the affine group ${\cal T}$-scheme $\pi ({\cal T})/ H_{\cal T}({\cal T}').$ The affine group $\pi ({\cal T})$-scheme 
$ H_{\cal T}({\cal T}')$ characterizes the Tannakian subcategory ${\cal T}'$ 
modulo $\otimes$-equivalence. In fact we have the following result 
(theorem 1.7):
 {\it there is bijection between the $\otimes$-equivalence classes of Tannakian subcategories of ${\cal T}$ and the normal affine group sub-${\cal T}$-schemes of $\pi ({\cal T})$, which associates
\par - to each Tannakian subcategory ${\cal T}' $ of 
$ {\cal T}$, the kernel $H_{\cal T}( {\cal T}')$ of the morphism of $ {\cal T}$-schemes $I: \pi( {\cal T}) \longrightarrow  i \pi( {\cal T}')$
corresponding to the inclusion 
$i:{\cal T}'\longrightarrow {\cal T};$ 
\par - to each  normal affine group sub-${\cal T}$-schemes $H$
 of $ \pi( {\cal T})$, the Tannakian subcategory 
${\cal T}(H)$ of ${\cal T}$, whose fundamental group $\pi({\cal T}(H))$ 
is  the affine group ${\cal T}$-scheme $\pi( {\cal T})/  H$}.

\pn Hence we obtain a clear dictionary between Tannakian subcategories of $\cal T$
and normal affine group sub-$\cal T$-schemes of the fundamental group $\pi({\cal T})$ of $\cal T$. Pierre Deligne has pointed out to the author that we can see this bijection as a ``reformulation'' of [11] II 4.3.2 b) and g). The proof of theorem 1.7 is based on [7] theorem 8.17.

\par We apply this dictionary to the Tannakian category ${\cal T}_1(k)$ generated by motives of
 niveau $\leq 1$ defined over $k$ (in an appropriate category of mixed realizations). 
We want to precise that in this article we restrict to motives of niveau $\leq 1$
because we are interested in motivic (and hence geometric) results and until now we know concretely only motives of niveau $\leq 1$.
The fundamental group of ${\cal T}_1(k)$ is called the {\it motivic Galois group 
$\G ({\cal T}_1(k))$ of motives of niveau $\leq 1$.}

\par For each fibre functor $\omega$ of ${\cal T}_1(k)$ 
over a $k$-scheme $S$, 
$\omega\G ({\cal T}_1(k))$ is the affine group $S$-scheme ${\underline {\rm Aut}}^{\otimes}_S(\omega)$ which represents the functor
which associates to each $S$-scheme
$T$, $u: T \longrightarrow S$, the group of automorphisms of
 $\otimes$-functors of the functor $u^* \omega .$ In particular, for each 
 embedding $\sigma:k \longrightarrow \CC$, the fibre functor $\omega_{\sigma}$
``Hodge realization'' furnishes the algebraic group $\QQ\,$-scheme 
$$\omega_{\sigma}  \G({\cal T}_1(k)) =  \spec \big(
\omega_{\sigma } (\Lambda) \big)= {\underline {\rm Aut}}^{\otimes}_{\QQ}
(\omega_{\sigma })\leqno(0.2)$$

\pn which is the {\it Hodge realization of the motivic Galois group of ${\cal T}_1(k)$.}

\par The weight filtration $W_*$ of the motives of niveau $\leq 1$ induces 
an increasing filtration $W_*$ of 3 steps in the 
affine group ${\cal T}_1(k)$-scheme $\G({\cal T}_1(k)).$ 
 Each of these 3 steps can be defined as intersection of normal affine group sub-${\cal T}_1(k)$-schemes of 
$\G({\cal T}_1(k))$, which correspond to determined Tannakian sub-categories of ${\cal T}_1(k)$.
For the Tannakian subcategories ${\Gr}^W_0 {\cal T}_1(k),{\Gr}^W_{-2} {\cal T}_1(k), 
{\Gr}^W_* {\cal T}_1(k)$ and $W_{-1}{\cal T}_1(k)$ 
of ${\cal T}_1(k)$, we compute the corresponding normal
affine group sub-${\cal T}_1(k)$-schemes which will furnish us, according to (0.1),  four exact short sequences of affine group sub-${\cal T}_1(k)$-schemes of 
$\G({\cal T}_1(k))$ (theorem 3.4). These short exact sequences are correlated one to each other by inclusions and projections, and they involve also the filtration $W_*$ of $\G({\cal T}_1(k))$.
One of these short exact sequences is 

$$ 0 \longrightarrow {\rm Res}_{\ok/k} \G({\cal T}_1(\ok)) \longrightarrow 
\G({\cal T}_{1}(k)) {\buildrel \pi \over  \longrightarrow } \gal  \longrightarrow 0.\leqno(0.3)$$

\pn If $\tau$ is an element of $\gal$, we have that $\pi^{-1}(\tau)$ is 
${\underline {\rm Hom}}^{\otimes}({\rm Id},\tau \circ {\rm Id} ),$ where 
$ {\rm Id}$ and $ \tau \circ {\rm Id}$ have to be regarded as functors on ${\cal T}_{1}(\ok)$ (corollary 3.6).
 By (0.2), the short exact sequence (0.3) is the motivic version of 
the short exact sequence 
of $\QQ \, $-algebraic groups

$$ 0 \longrightarrow {\underline {\rm Aut}}^{\otimes}_\QQ(\omega_{ {\overline \sigma } \, | {\cal T}_1(\ok)}) \longrightarrow 
{\underline {\rm Aut}}^{\otimes}_\QQ (\omega_{\sigma })  \longrightarrow  \gal  \longrightarrow 0
\leqno(0.4)$$

\pn where ${\overline \sigma }: \ok \longrightarrow \CC$ is the embedding of $\ok$ in $\CC$ which extends $\sigma: k \longrightarrow \CC$. This last sequence is the restriction 
to motives of niveau $\leq 1$ of the sequence 
 found by P. Deligne and U. Jannsen in [5] II 6.23 and [8] 4.7 respectively.
Moreover the equality
 $\pi^{-1}(\tau) = {\underline {\rm Hom}}^{\otimes}({\rm Id},\tau \circ {\rm Id} )$
 is the motivic version of the one found by P. Deligne and U. Jannsen (loc. cit.).

\par Let $M$ be a 1-motive defined over $k$. The motivic Galois group $\G(M)$ of $M$ is the fundamental group of the Tannakian subcategory $\langle M \rangle^\otimes$ of ${\cal T}_1(k)$ generated by $M$. In [2], we compute the unipotent radical $W_{-1}(\L \G (M))$ of the Lie algebra of $\G(M)$: it is the semi-abelian variety defined by the adjoint action of the graduated  
${\Gr}_*^W(W_{-1} \L \G (M))$ on itself. This result allows us to compute the derived group of the unipotent radical $W_{-1}(\L \G (M))$ (proposition 4.4).  
Applying then theorem 1.7, we construct explicitly the biggest Tannakian subcategory of $\langle M \rangle^\otimes$ which has a commutative motivic Galois group (theorem 4.9).

\par The ``Hodge realization'' of the motivic Galois group of  motives was partially studied by 
P. Deligne and U. Jannsen in [5] II \S 6 and [8] 4.6 respectively. 
They don't restrict to motives of niveau $\leq 1$ and hence we find the motivic version 
of the restriction to motives of niveau $\leq 1$ of some of theirs results (corollary 3.6).

\par The dictionary between Tannakian subcategories of a Tannakian category $\cal T$
 and normal affine group sub-$\cal T$-schemes of the fundamental group $\pi({\cal T})$ of $\cal T$ has applications also in a ``non-motivic'' context. It can be apply whenever one has objects generating a Tannakian category.

\par In the first section, we recall the definition and some properties of the fundamental group 
$\pi( {\cal T})$ of a Tannakian category ${\cal T}$ and we prove the bijection between 
the $\otimes$-equivalence classes of Tannakian subcategories of ${\cal T}$ 
and normal affine sub-${\cal T}$-schemes of $\pi( {\cal T})$.  
\pn In the section 2, we apply the definitions of the first section to the Tannakian category
 ${\cal T}_1(k)$
generated by motives of niveau $\leq 1$ defined over $k$ 
(in an appropriate category of mixed realizations): the motivic Galois
 group of a motive is the fundamental group of the Tannakian category generated by this motive.
 We end this section giving several examples of motivic Galois groups. 
\pn In the third section, we apply the dictionary between Tannakian sub-categories
 and normal affine group sub-${\cal T}_1(k)$-schemes to some Tannakian subcategories of ${\cal T}_1(k).$
 In particular, we get the motivic version of [5] II 6.23 (a), (c) and [8] 4.7 (c), (e).  
\pn In section 4, we compute the biggest Tannakian subcategory of the Tannakian category generated by a 1-motive, which has a commutative 
motivic Galois group.

\vskip 0.5 true cm
\par {\it Acknowledgements:} I want to express my gratitude to Pierre Deligne 
for his precisions about the notion of Tannakian subcategory and about theorem 1.7. I thank Uwe Jannsen for the various discussions we have had on 1-motives. I thank also the department 
of mathematics of Regensburg for the very nice conditions of work.

\vskip 0.5 true cm
\pn In this paper $k$ is a field of characteristic 0 embeddable in $\CC$ and $\ok$ its algebraic closure.
\vskip 0.5 true cm

\pn {\bf 1. Motivic Galois theory.}
\vskip 0.5 true cm

\pn 1.1. Let $\cal T$ be a Tannakian category over $k$, i.e. a
tensor category over $k$ with a fibre functor over a nonempty $k$-scheme.
A {\bf Tannakian subcategory
of $\cal T$} is a strictly full subcategory ${\cal T}'$ of $\cal T$ which is 
closed under the formation of subquotients, direct sums, tensor products and duals,
and which is endowed with the restriction to ${\cal T}'$ 
of the fibre functor of $\cal T$.  

\pn The tensor product of $\cal T$ extends to a tensor product in the category 
${\rm Ind}{\cal T}$ of Ind-objects of $\cal T$.
 A {\bf commutative ring} in ${\rm Ind}{\cal T}$ is an object $A$ of
 ${\rm Ind}{\cal T}$ together with a commutative associative product 
$A \otimes A \longrightarrow A$ admitting an identity
 $1_{\cal T} \longrightarrow A,$
where $1_{\cal T}$ is the unit object of $\cal T$.
 In ${\rm Ind}{\cal T}$ we can define in the usual way the notion of morphisms of commutative rings,
the notion of $A$-module (an object of 
${\rm Ind}{\cal T}$ endowed with a structure of $A$-module), ...
The {\bf category of affine $\cal T$-schemes} is the opposite
of the category of commutative rings in ${\rm Ind}{\cal T}$.
We denote ${\rm Sp}(A)$ the affine $\cal T$-scheme defined by the ring $A$.
 The final object of this category is the affine $\cal T$-scheme  
${\rm Sp}(1_{\cal T})$ defined by the ring $1_{\cal T}$.
 A module over ${\rm Sp}(A)$ is an $A$-module and for each 
morphism ${\rm Sp}(B) \longrightarrow {\rm Sp}(A)$, the functor 
$M \longmapsto B \otimes_A M$ is called ``the inverse image over 
${\rm Sp}(B).$''
An {\bf affine group $\cal T$-scheme} is a group object in the category 
of affine $\cal T$-schemes, i.e. ${\rm Sp}(A)$ with $A$ endowed with a structure of Hopf algebra. An {\bf action} of an affine group $\cal T$-scheme ${\rm Sp}(A)$ on an object $X$ of $\cal T$ is a morphism 
$X \longrightarrow X \otimes A$ satisfying the usual axioms for a $A$-comodule.
The {\bf Lie algebra} of an affine group $\cal T$-scheme 
is a  pro-object $\rm L$ of ${\cal T}$ endowed with a Lie algebra structure.
 
\vskip 0.5 true cm

\pn 1.2. REMARKS: 
\par (1) Since ${\rm End}(1_{\cal T})=k,$ the Tannakian subcategory of $\cal T$ consisting
of direct sum of copies of $1_{\cal T}$ is equivalent to the Tannakian category 
of finite-dimensional $k$-vector spaces. In particular, {\it every Tannakian category over $k$ contains, as Tannakian subcategory, the category of finite-dimensional $k$-vector spaces ${\rm Vec}(k)$}.
Considering Ind-objects, we have that 
{\it each affine $k$-scheme defines an affine $\cal T$-scheme}. In particular,
${\rm Spec} (k)$ corresponds to ${\rm Sp}(1_{\cal T})$ (cf. [6] 5.6).
\par (2) Let $\cal T$ be the Tannakian category of representations of an affine group scheme $G$ over $k$: ${\cal T}= {\rm Rep}_k(G)$. In this case, 
affine $\cal T$-schemes are affine $k$-schemes endowed with an action of $G$.
 The inclusion of affine $k$-schemes in the category of
affine $\cal T$-schemes described in 1.2 (1), is realized adding the trivial action of $G$ 
(cf. [6] 5.8).
\vskip 0.5 true cm

\pn 1.3. The {\bf fundamental group $\pi({\cal T})$} of a Tannakian category
${\cal T}$ is the affine group $\cal T$-scheme
 ${\rm Sp}(\Lambda),$ where
 $ \Lambda$ is the Hopf algebra of ${\rm Ind}{\cal T}$ satisfying the following universal property: for each object $X$ of ${\cal T}$,
there exists a morphism
$$\lambda_X: X^{\du} \otimes X \longrightarrow \Lambda \leqno(1.3.1)$$
\pn functorial in $X$, i.e.
for each $f:X \longrightarrow Y$ in ${\cal T}$ the diagram

$$\matrix{Y^{\du} \otimes X&{\buildrel f^t \otimes 1 \over \longrightarrow}&
X^{\du} \otimes X\cr
     {\scriptstyle 1 \otimes f}\downarrow~~~~ && ~~~~~ \downarrow
 {\scriptstyle \lambda_X}  \cr
   Y^{\du} \otimes Y& {\buildrel \lambda_Y \over \longrightarrow}&\Lambda \cr}$$

\pn is commutative. The universal property of $\Lambda$ reads:
for each Ind-object $U$ of ${\cal T}$ the application

 $$\eqalign{{\Hom}(\Lambda, U)& \longrightarrow \{ u_X:
  X^{\du} \otimes X \longrightarrow U, ~~ {\rm {functorial~ in~}} X \} \cr
f & \longmapsto  f  \circ \lambda_X \cr }\leqno(1.3.2)$$
 
\pn is a bijection. The existence of the fundamental group 
$\pi ({\cal T}) $ is proved in [7] 8.4, 8.10, 8.11 (iii). 
\par The morphisms (1.3.1), which  can be rewritten on the form

 $$ X \longrightarrow  X \otimes \Lambda,\leqno(1.3.3)$$

\pn define {\it an action of the fundamental group $\pi ({\cal T}) $ on each object of ${\cal T}$.} 
 In parti\-cu\-lar, the morphism 
$ \Lambda \longrightarrow  \Lambda \otimes \Lambda$
represents the action of $\pi ({\cal T}) $ on itself by inner automorphisms
(cf. [6] 6.1).
\vskip 0.5 true cm

\pn 1.4. EXAMPLES:
\par (1) Let ${\cal T}={\rm Vec}(k)$ be the Tannakian category of
 finite dimensional $k$-vector spaces. From the main theorem of
 Tannakian category, we know that ${\rm Vec}(k)$ is equivalent to
 the category of finite-dimensional $k$-representations of $\spec(k)$. 
In this case, affine ${\cal T}$-schemes are affine $k$-schemes
and $\pi ({\rm Vec}(k)) $ is $\spec(k)$.
\par (2)  Let $\cal T$ be the Tannakian category of $k$-representations
 of an affine group scheme $G$ over $k$: ${\cal T}= {\rm Rep}_k (G)$.
 From [6] 6.3, the fundamental group $\pi ({\cal T}) $ of ${\cal T}$
 is the affine group $k$-scheme $G$
endowed with its action on itself by inner automorphisms.
\vskip 0.5 true cm

\pn 1.5. From [7] 6.4, to any exact and $k$-linear $\otimes$-functor
$u: {\cal T}_1 \longrightarrow  {\cal T}_2$ between Tannakian categories over $k,$ 
corresponds a morphism of affine group $ {\cal T}_2$-schemes

$$U: \pi( {\cal T}_2) \longrightarrow u \pi( {\cal T}_1). \leqno(1.5.1)$$

\pn For each object $X_1$ of ${\cal T}_1$,
the action (1.3.3) of $\pi( {\cal T}_1)$ on $X_1$ induces an action of
$ u \pi( {\cal T}_1)$ on $u(X_1)$. Via (1.5.1) this last action 
induces the action (1.3.3) of $\pi( {\cal T}_2)$ on the object 
$u(X_1)$ of ${\cal T}_2$.

\par By the theorem [11] II 4.3.2 (g) we have the following dictionary 
between the functor $u$ and the morphism $U$:
\par (1) $U$ is faithfully flat (i.e. flat and surjective) if and only if
$u$ is fully faithful and every subobject of $u(X_1)$ for $X_1$ an object of $ {\cal T}_1$, is isomorphic to the image of a subobject of $X_1$.
\par (2) $U$ is a closed immersion if and only if every object of 
 $ {\cal T}_2$ is isomorphic to a subquotient of an object of the form 
$u(X_1)$, for $X_1$ an object of $ {\cal T}_1$.
\vskip 0.5 true cm

\pn 1.6. We can now state the dictionary between Tannakian subcategories 
of ${\cal T}$ and normal affine group sub-${\cal T}$-schemes 
of the fundamental group $\pi({\cal T})$ of ${\cal T}$:
\vskip 0.5 true cm

\pn 1.7. {\bf Theorem}
\pn {\it Let ${\cal T}$ be a Tannakian category over $k$, with fundamental group 
$ \pi( {\cal T})$.
There is a bijection between the $\otimes$-equivalence classes of 
Tannakian subcategories 
of ${\cal T}$ and the normal affine group sub-${\cal T}$-schemes
 of $ \pi( {\cal T}):$

$$\eqalign{ 
\pmatrix{\otimes-{\rm equiv.~ classes~ of~}\cr
{\rm Tannakian~ subcat.~  of~} {\cal T} \cr } 
& \longrightarrow 
\pmatrix{{\rm normal~ affine~ group~}\cr
{\rm sub}-{\cal T}-{\rm schemes~ of}~  \pi( {\cal T}) \cr }  \cr
 {\cal T}' &~ \longrightarrow ~  H_{\cal T}  ( {\cal T}')=\ker\big( \pi( {\cal T}){\buildrel I \over \longrightarrow} i \pi( {\cal T}')\big)\cr
{\cal T}(H)~{\rm s.t.~}\pi({\cal T}(H)) \cong \pi({\cal T})/H &~
 \longleftarrow ~  H \cr}$$

\pn which associates 

\par - to each Tannakian subcategory ${\cal T}' $ of 
$ {\cal T}$, the kernel $H_{\cal T}( {\cal T}')$ of the morphism of affine group $ {\cal T}$-schemes
$I: \pi( {\cal T}) \longrightarrow  i \pi( {\cal T}')$
corresponding to the inclusion functor 
$i:{\cal T}'\longrightarrow {\cal T}.$ In particular,
 we can identify the fundamental group $\pi( {\cal T}')$ 
of $ {\cal T}'$ with the affine group ${\cal T}$-scheme
$\pi( {\cal T})/  H_{\cal T}  ( {\cal T}').$

\par - to each  normal affine group sub-${\cal T}$-scheme $H$
 of $ \pi( {\cal T})$, the Tannakian subcategory 
${\cal T}(H)$ of ${\cal T}$ whose fundamental group $\pi({\cal T}(H))$ 
is  the affine group ${\cal T}$-scheme $\pi( {\cal T})/  H.$
 Explicitly, ${\cal T}(H)$ is the Tannakian subcategory of ${\cal T}$
generated by the objects 
of ${\cal T}$ on which the action (1.3.3) of $ \pi( {\cal T})$ induces
 a trivial action of $H.$}
\vskip 0.5 true cm

\pn PROOF: 
 Let $i: {\cal T}' \longrightarrow  {\cal T}$
be the inclusion functor of a Tannakian sub-category ${\cal T}'$ of ${\cal T}$.
 From 1.5 to this functor corresponds 
 a faithfully flat morphism of affine group $ {\cal T}$-schemes
$I: \pi( {\cal T}) \longrightarrow  i \pi( {\cal T}').$
Denote by 

$$H_{\cal T}  ( {\cal T}')=\ker\big( \pi( {\cal T})
{\buildrel I \over \longrightarrow} i \pi( {\cal T}')\big)$$ 

\pn the affine group $ {\cal T}$-scheme which is the kernel of 
the morphism $I$. In particular, we have the exact sequence of affine group $ {\cal T}$-schemes

$$0 \longrightarrow H_{\cal T}( {\cal T}') \longrightarrow \pi( {\cal T})
\longrightarrow i \pi( {\cal T}') \longrightarrow 0.$$

\pn According to theorem 8.17 of [7], the inclusion functor
$i:{\cal T}'\longrightarrow {\cal T}$ identifies 
 ${\cal T}'$ with the Tannakian subcategory of objects of ${\cal T}$ on which
the action (1.3.3) of $ \pi( {\cal T})$ induces a trivial action of
$H_{\cal T}( {\cal T}')$. In particular, we get that

$$\pi( {\cal T}') \cong \pi( {\cal T})/  H_{\cal T}  ( {\cal T}').
\leqno(1.7.1)$$

\pn We now check the injectivity: If ${\cal T}_1$ and ${\cal T}_2$ are two 
Tannakian sub-categories of ${\cal T}$ such that $H_{\cal T} ( {\cal T}_1)=
H_{\cal T}  ( {\cal T}_2)$, then by (1.7.1)
they have also the same fundamental group. But this means that there is a 
$\otimes$-equivalence of categories between ${\cal T}_1$ and ${\cal T}_2$.
\pn The surjectivity is trivial.
\vskip 0.5 true cm

\pn REMARK: The existence of quotient affine group ${\cal T}$-schemes is assured by 
5.14 (ii) [6].
\vskip 0.5 true cm

\pn {\bf 1.8. Lemma}
\pn {\it Let ${\cal T}$ be a Tannakian category over $k$, with fundamental group $\pi ({\cal T})$.
\par (i) If ${\cal T}_1,{\cal T}_2$ are two Tannakian subcategories of $ {\cal T}$ such that ${\cal T}_1 \subseteq {\cal T}_2$, 
then $H_{\cal T}({\cal T}_1) \supseteq H_{\cal T}({\cal T}_2).$
\par (ii) If $H_1,H_2$ are two normal subgroups of $\pi ({\cal T})$ such that $H_1 \subseteq H_2$, 
then ${\cal T}(H_1) \supseteq {\cal T}(H_2).$ In particular we have the projection
$$\pi( {\cal T}(H_1))= \pi ({\cal T})/H_1 \longrightarrow \pi( {\cal T}(H_2))= \pi ({\cal T})/H_2.$$}

\pn 1.9. EXAMPLE: Let $ {\cal T}$ be a Tannakian category over $k$.
The inclusion $i: {\rm Vec}(k) \hookrightarrow {\cal T}$ 
of the Tannakian category $ {\rm Vec}(k)$ of finite-dimensional $k$-vector spaces in $ {\cal T}$, defines the faithfully flat morphism of affine group $ {\cal T}$-schemes
$$I: \pi( {\cal T}) \longrightarrow i \spec (k).$$
\pn whose kernel is the whole $ {\cal T}$-scheme $\pi( {\cal T})$: $ H_{\cal T}  ({\rm Vec}(k) )=\pi( {\cal T})$.
Hence the functor $i$ identifies the category $ {\rm Vec}(k)$  
 with the Tannakian subcategory of objects of ${\cal T}$ on which
the action (1.3.3) of $ \pi( {\cal T})$ is trivial (cf. [6] 6.7 (i)). 
\vskip 0.5 true cm

\pn 1.10. Let $\omega$ be a fibre functor of 
the Tannakian category ${\cal T}$ over a 
 $k$-scheme $S$, namely an exact $k$-linear $\otimes$-functor from
${\cal T}$ to the category of quasi-coherent sheaves over $S$.
 It defines a $\otimes$-functor, denoted again $\omega$, 
from ${\rm Ind}{\cal T}$ to the category of quasi-coherent sheaves over $S$.
If $\pi ({\cal T})={\rm Sp}(\Lambda)$ we define

 $$\omega(\pi ({\cal T}))= {\rm Spec} (\omega (\Lambda)).\leqno(1.10.1)$$

\pn According [7] (8.13.1), the spectrum  
${\rm Spec} (\omega (\Lambda))$ is the affine group $S$-scheme 
$ {\underline {\rm Aut}}^{\otimes}_S(\omega)$ which represents the functor
which associates to each $S$-scheme
$T$, $u: T \longrightarrow S$, the group of automorphisms of
 $\otimes$-functors of the functor 

$$\eqalign{ \omega_T: {\cal T}& \longrightarrow  \{ {\rm locally~ free~  scheaves~ of~ finite~ rank~ over~} T \} \cr 
 X & \longmapsto u^*\omega(X).\cr }$$

\pn From the formalism of [6] 5.11, we have the following dictionary:

\par - to give the affine group ${\cal T}$-scheme
$\pi ({\cal T})={\rm Sp}(\Lambda)$ 
is the same thing as to give, for each fibre functor
 $\omega$ over a $k$-scheme $S$, the affine group  
$S$-scheme $ {\underline {\rm Aut}}^{\otimes}_S(\omega)$, in a functorial way with respect to $\omega$ and in a compatible way with respect to the base changes 
$S' \longrightarrow S.$

\par - let $u: {\cal T}_1 \longrightarrow  {\cal T}_2$ be a 
$k$-linear $\otimes$-functor between Tannakian categories over $k.$ 
To give the corresponding morphism 
$U: \pi( {\cal T}_2) \longrightarrow u \pi( {\cal T}_1)$
 of affine group $ {\cal T}_2$-schemes,
is the same thing as to give, for each fibre functor $\omega$ of $ {\cal T}_2$
 over a $k$-scheme $S$, a morphism of affine group $S$-schemes
$ {\underline {\rm Aut}}^{\otimes}_S(\omega) \longrightarrow 
{\underline {\rm Aut}}^{\otimes}_S(\omega \circ u)$, in a functorial way with respect to $\omega$.
\vskip 0.5 true cm

\pn 1.11. {\bf Lemma}
\pn {\it  Let $u_1: {\cal T}_1 \longrightarrow {\cal T}_2$ and 
$ u_2: {\cal T}_2 \longrightarrow {\cal T}_3$ be two exact and $k$-linear $\otimes$-functors between Tannakian categories over $k$. Denote by 
$U_1: \pi({\cal T}_2) \longrightarrow u_1 \pi({\cal T}_1)$ and 
$U_2: \pi({\cal T}_3) \longrightarrow u_2 \pi({\cal T}_2)$
the morphisms of affine group ${\cal T}_2$-schemes and 
${\cal T}_3$-schemes defined respectively by $u_1$ and $u_2$.
Then the morphism of affine group ${\cal T}_3$-schemes corresponding to 
$u_2 \circ u_1$ is 

$$ U = u_2 U_1 \circ U_2 :\pi({\cal T}_3) \longrightarrow u_2 \pi({\cal T}_2) \longrightarrow u_2 u_1 \pi({\cal T}_1).  $$

\pn Moreover,
\par (i) if $u_2 \circ u_1 \equiv 1_{{\cal T}_3}$ then $U: \pi({\cal T}_3)
\longrightarrow {\sp}(1_{{\cal T}_3}),$
\par (ii) if ${\cal T}_1={\cal T}_2={\cal T}_3$ and $u_2 \circ u_1=id$, then 
$U=id$.}
\vskip 0.5 true cm

\pn PROOF: The morphism of affine group ${\cal T}_2$-schemes 
$U_1: \pi({\cal T}_2) \longrightarrow u_1 \pi({\cal T}_1)$ furnishes 
a morphism of affine group ${\cal T}_3$-schemes

$$u_2 U_1:u_2 \pi({\cal T}_2) \longrightarrow u_2 u_1 \pi({\cal T}_2)$$

\pn which corresponds to the following systeme of morphisms: 
for each fibre functor $\omega$  of $ {\cal T}_3$
 over a $k$-scheme $S$, we have a morphism of affine group $S$-schemes

$$ {\underline {\rm Aut}}^{\otimes}_S(\omega \circ u_2) \longrightarrow 
{\underline {\rm Aut}}^{\otimes}_S\big((\omega \circ u_2) \circ u_1\big).
\leqno(1.11.1)$$
 
\pn Denote by $U:\pi({\cal T}_3) \longrightarrow u_2 u_1 \pi({\cal T}_1)$ the
morphism of affine group ${\cal T}_3$-schemes corresponding to the functor
$ u_2 \circ u_1:{\cal T}_1 \longrightarrow {\cal T}_3$. To have the morphism 
$U$ (resp. $U_2$) of ${\cal T}_3$-schemes is the same thing as to have, 
for each fibre functor $\omega$  of $ {\cal T}_3$
 over a $k$-scheme $S$, a morphism of affine group $S$-schemes

$$\eqalign{
 {\underline {\rm Aut}}^{\otimes}_S(\omega )& \longrightarrow 
{\underline {\rm Aut}}^{\otimes}_S\big(\omega \circ (u_2\circ u_1)\big)\cr
({\rm resp.}~~~
  {\underline {\rm Aut}}^{\otimes}_S(\omega )& \longrightarrow 
{\underline {\rm Aut}}^{\otimes}_S\big(\omega \circ u_2 \big)~~)\cr}\leqno(1.11.2)$$
 
\pn Hence, according to (1.11.1) we observe that $U = u_2 U_1 \circ U_2.$
\pn The remaining assertions are clear from (1.11.2): in particular, if 
$u_2 \circ u_1 \equiv 1_{{\cal T}_3},$ we have that ${\underline {\rm Aut}}^{\otimes}_S\big(\omega_{\vert \langle 1_{{\cal T}_3} \rangle^\otimes}
\big)= {\spec} (k)$ for each fibre functor $\omega$  of $ {\cal T}_3$
 over a $k$-scheme $S$.
\vskip 0.5 true cm

\pn {\bf 2. Some motivic Galois groups.}
\vskip 0.5 true cm

\pn 2.1. Let $MR(k)$ be the category of mixed realizations
 (for absolute Hodge cycles) over $k$ defined by U. Jannsen in [8] I 2.1. 
 It is a neutral Tannakian category over $\QQ$, 
i.e. a tensor category over $\QQ$ with a fibre functor over $\spec(\QQ).$ 
Each embedding $\sigma:k \longrightarrow \CC$ gives a fibre functor $\omega_\sigma$ of $MR(k)$ over $\spec(\QQ)$, called ``the Hodge realization''. 
If $k=\QQ$, there is another fibre functor $\omega_{\rm dR}$ of $MR(k)$ over $\spec(\QQ)$, called ``the de Rham realization''.
\par  {\bf The category of motives over $k$} is the Tannakian
subcategory of $MR(k)$ generated (by means of $\oplus$, $\otimes$, dual,
subquotients) by those realizations which come from geometry (cf. [6] 1.11).

\pn A motive over $k$ {\bf with integral coefficients} is a motive $M$
 whose mixed realization is endowed with an integral structure, 
i.e. for each embedding $\sigma:k \rightarrow \CC$, the Hodge realization $
{\rm T}_\sigma(M)$ of $M$ contains a $\ZZ$-lattice ${\rm L}_\sigma$
and for each prime number $\ell$, the $\ell$-adic realization  $
{\rm T}_{\ell} (M)$ of $M$ contains a $\gal$-invariant $\ZZ_\ell$-lattice
 ${\rm L}_\ell ,$
 such that $I_{{\overline \sigma}, {\ell}}({\rm L}_\sigma)={\rm L}_\ell$
 for each $\ell$ and ${\overline \sigma}$ (cf. [6] 1.23). 1-Motives over
 $k$ are motives with integral coefficients.

\par In this paper we are interested in two kinds of motives over $k$: 
Artin motives and 1-motives.
\vskip 0.5 true cm

\pn 2.2. By [5] 6.17 we have a fully faithful functor
from the category of 0-dimensio\-nal varieties over $k$ to the category 
$MR(k)$ of mixed realizations over $k$

$$\eqalign{ \{ {\rm 0-dim.~varieties}/k \} & \longrightarrow MR(k)\cr
 X & \longmapsto ({\rm T}_\sigma(X),{\rm T}_{\rm dR}(X),{\rm T}_{\ell}(X),
I_{\sigma, {\rm dR}}, I_{{\overline \sigma}, {\ell}} )_{{\sigma:k \rightarrow \CC, ~ {\overline \sigma}:\ok \rightarrow \CC} \atop
{{\ell}~{\rm prime ~ number}}}\cr}$$

\pn which associates to each 0-dimensional varietie $X$ its Hodge, de Rham and
 $\ell$-adic realizations and the comparison isomorphisms
 $I_{\sigma, {\rm dR}}:{\rm T}_\sigma(X) \otimes_\QQ \CC \longrightarrow 
{\rm T}_{\rm dR}(X)\otimes_k \CC $ and $I_{{\overline \sigma},
 {\ell}}:{\rm T}_\sigma(X) \otimes_\QQ 
\QQ_\ell \longrightarrow {\rm T}_{\ell}(X)$.
Hence we can identify 0-dimensional varieties with the mixed realizations they define.

\par The {\bf Tannakian category of Artin motives ${\cal T}_{0}(k)$ 
over $k$} is the Tannakian subcategory of $MR(k)$
generated by 0-dimensional varieties over $k$, i.e. by
mixed realizations of 0-dimensional varietes.
{\it ${\cal T}_{0}(k)$ is a neutral Tannakian category
 over $\QQ$ with fibre functors $\{ \omega_\sigma \}_{\sigma:k \longrightarrow \CC}$
``the Hodge realizations''.}
The unit object of ${\cal T}_{0}(k)$ is $\spec(k)$.
\pn Through the Hodge realization,
 ${\cal T}_{0}(k)$ is equivalent to the category 
of finite-dimen\-sio\-nal $\QQ \,$-representations of $\gal$, i.e.

$$\eqalign{ {\cal T}_{0}(k)& ~~ \cong ~~ {\rm Rep}_{\QQ}(\gal)\cr
X &~ \longmapsto ~ \QQ^{X(\ok)}.\cr}\leqno(2.2.1)$$

\pn Here and in the following we regard $\gal$ as a constant, pro-finite affine group scheme over $\QQ$.
 It is clear that ${\cal T}_{0}(\ok)$ is equivalent to 
the Tannakian category of 
finite-dimensional $\QQ \,$-vector spaces : in fact, as affine group $\QQ\,$-scheme ${\rm Gal}(\ok / \ok)$ is 
$\spec (\QQ)$, and so

$${\cal T}_{0}(\ok) \cong {\rm Rep}_{\QQ}( \spec (\QQ)) \cong {\rm Vec}(\QQ).
\leqno(2.2.2)$$

\par Artin motives are pure motives of weight 0 : the weight 
 filtration $W_*$ on an Artin motive $X$ is ${\W}_{i}(X) = X$ for each $i \geq 0$ and ${\W}_{j}(X) = 0$ for each $j \leq -1$. If we denote ${\rm Gr}_{n}^{{\W}}\, =\,
{\W}_{n} / {\W}_{n-1},$ we have ${\rm Gr}_{0}^{{\W}}(X)=X$ and 
$ {\rm Gr}_{i}^{{\W}}(X)=0$ for each $i \not= 0$.
\vskip 0.5 true cm

\pn 2.3. A {\bf 1-motive} $M$ over $k$ consists of
\vskip 0.2 true cm

\par (a) a group scheme $X$ over $k$, which is locally for the \'etale
topology, a constant group scheme defined by a finitely generated free
$\ZZ$-module,
\vskip 0.2 true cm

\par (b) a semi-abelian variety $G$ defined over $k,$ i.e. 
an extention of an abelian variety $A$ by a torus $Y(1),$
which cocharacter group $Y$,
\vskip 0.2 true cm

\par (c) a morphism $u:X \longrightarrow G$ of group schemes over $k$.
\vskip 0.2 true cm

\par We have to think of $X$ as a character group of a torus
 defined over $k$, i.e. as a finitely generated 
$\gal$-module. We identify $X(\ok)$ with a free 
 $\ZZ$-module finitely generated.
The morphism $u:X \longrightarrow G$ is equivalent to a 
$\gal$-equivariant homomorphism $u:  X(\ok) \longrightarrow G(\ok).$

\par An {\bf isogeny } between two 1-motives
$M_{1}=[X_{1} {\buildrel u_{1} \over \longrightarrow} G_{1}]$ and 
$M_{2}=[X_{2} {\buildrel u_{2} \over \longrightarrow} G_{2}]$ is a morphism 
of 1-motives (i.e. 
 a morphism of complexes of commutative group schemes)
such that 
 $f_{X}:X_{1} \longrightarrow X_{2}$ is injective with finite cokernel, and 
 $f_{G}:G_{1} \longrightarrow G_{2}$ is surjective with finite kernel.

\par 1-motives are mixed motives of niveau $\leq 1$: the weight filtration $W_*$ on $M
=[X {\buildrel u \over \longrightarrow} G] $ is 

$$ \eqalign{{\W}_{i}(M)&=M  ~~{\rm  for ~ each~} i \geq 0, \cr   
{\W}_{-1}(M)& =[0 \longrightarrow G],\cr
 {\W}_{-2}(M)&=[0 \longrightarrow  Y(1)],\cr 
  {\W}_{j}(M) & =0 ~~~{\rm  for ~ each~} j \leq -3. \cr}$$

\pn In particular, we have  ${\rm Gr}_{0}^{{\W}}(M)= 
[X {\buildrel  \over \longrightarrow} 0], {\rm Gr}_{-1}^{{\W}}(M)= 
[0 {\buildrel  \over \longrightarrow} A]$ and $ {\rm Gr}_{-2}^{{\W}}(M)= 
[0 {\buildrel  \over \longrightarrow}  Y(1)].$
\vskip 0.5 true cm

 \pn 2.4. According to [4] 10.1.3 (other [9] 4.2 (i) other for $k=\QQ$ [6] 2.3) we have a fully faithful functor
from the category  of 1-motives over $k$ to the category $MR(k)$ of mixed 
realizations over $k$

$$\eqalign{ \{ {\rm 1-motives} /k \}& \longrightarrow MR(k)\cr 
 M & \longmapsto ({\rm T}_\sigma(M),{\rm T}_{\rm dR}(M),{\rm T}_{\ell}(M),
I_{\sigma, {\rm dR}}, I_{{\overline \sigma}, {\ell}} )_{{\sigma:k \rightarrow \CC, ~ {\overline \sigma}:\ok \rightarrow \CC} \atop
{{\ell}~{\rm prime ~ number}}}\cr}$$

\pn which associates to each 1-motive $M$ its Hodge, de Rham and
 $\ell$-adic realizations and the comparison isomorphisms
 $I_{\sigma, {\rm dR}}:{\rm T}_\sigma(M) \otimes_\QQ \CC \longrightarrow 
{\rm T}_{\rm dR}(M)\otimes_k \CC $ and $I_{{\overline \sigma}, {\ell}}:{\rm T}_\sigma(M) \otimes_\QQ 
\QQ_\ell \longrightarrow {\rm T}_{\ell}(M)$.
 Therefore we can identify 1-motives with the mixed realizations they define.

\par The {\bf Tannakian category ${\cal T}_{ 1}(k)$ of 1-motives over $k$} is
the Tannakian subcategory of
$MR(k)$ generated by 1-motives, i.e. by mixed realizations of 1-motives:
{\it  ${\cal T}_{1}(k)$
is a neutral Tannakian category over $\QQ$ with fibre functors $\{ \omega_\sigma
\}_{\sigma : k \longrightarrow \CC}$
``the Hodge rea\-li\-za\-tions''.}
The unit object of ${\cal T}_{ 1}(k)$ is the 1-motive 
$\ZZ(0)=[\ZZ  \longrightarrow 0]$. 
For each object $M$ of ${\cal T}_1(k)$, we denote by 
$M^{\du}={\underline {\Hom}}(M,\ZZ(0))$ its dual. 
The Cartier dual of an object $M$ of ${\cal T}_{ 1}(k)$ is the object 
$$M^* =M^{\du} \otimes \ZZ(1)\leqno(2.4.1)$$
\pn of ${\cal T}_{ 1}(k).$ If $M$ is a 1-motive, then  
$M^*$ is again a 1-motive.

\par We will denote by $W_{-1}{\cal T}_{1}(k)$ (resp. ${\Gr}^W_0 {\cal T}_{1}(k)$, ...)
the Tannakian subcategory of ${\cal T}_{1}(k)$ generated by the mixed realizations of all $W_{-1}M$
(resp. ${\Gr}^W_0 M$, ...) with $M$ a 1-motive.
\vskip 0.5 true cm

\pn 2.5. If a Tannakian category ${\cal T}$ is generated by motives, 
the fundamental group  $\pi ({\cal T}) $ is called 
the {\bf motivic Galois group $\G( {\cal T})$ of ${\cal T}:$} hence, it is the affine group 
${\cal T}$-scheme ${\sp}(\Lambda)$, where $\Lambda$ is the Hopf algebra of 
${\rm Ind}{\cal T}$ satisfying the universal property (1.3.1).
If the Tannakian category is generated by only one motive $M$,
i.e. ${\cal T}=\langle M \rangle^\otimes,$ $\pi ({\cal T}) $ 
is the {\bf motivic Galois group $\G(M)$ of $M$:} in this case,
 it is the affine group ${\cal T}$-scheme ${\sp}(\Lambda)$, where
 the Hopf algebra $\Lambda$ satisfying the universal property (1.3.1) 
 is an object of ${\cal T}$. (cf. [7] 6.12). 

\par Since the fundamental group depends only on the Tannakian 
category, the notion of motivic Galois group is stable under isogeny and duality.

\par After identifying motives with their mixed realizations and after
 taking the Tannakian category defined by those mixed realizations, 
we lose the integral structure 
(in fact the Tannakian categories we consider are $\QQ\,$-linear). This implies
 that the motivic Galois group, which depends only on the Tannakian category
 generated by motives, doesn't ``see'' the integral structure of motives. 
Hence, we will work with iso-motives. 
\vskip 0.5 true cm

\pn 2.6. EXAMPLES: 
\par (1) The motivic Galois group $\G(\ZZ(0))$ of the unit object $\ZZ(0)$ of
 ${\cal T}_1(k)$ is the affine group $\langle \ZZ(0) \rangle^\otimes$-scheme ${\rm Sp}(\ZZ(0))$.
For each fibre functor ``Hodge realization'' $\omega_\sigma$,
we have that
$ \omega_\sigma \big( \G(\ZZ(0))\big):={\rm Spec}\big(\omega_\sigma (\ZZ(0)) \big)=
{\rm Spec}( \QQ)$, which is the Mumford-Tate group of ${\rm T}_\sigma ( \ZZ(0))$.
\par (2) Let $\langle \ZZ(1) \rangle^\otimes$ be the neutral Tannakian
 category over $\QQ$ defined by the $k$-torus $\ZZ(1)$.
The motivic Galois group $\G (\ZZ(1))$ of the torus $\ZZ(1)$ 
 is the affine group $\langle \ZZ(1) \rangle^\otimes$-scheme ${\cal G}_{m}$ defined by the affine 
$\QQ\,$-scheme ${\GG}_{m / \QQ}$ (cf. remark 1.2 (1)).
For each fibre functor ``Hodge realization'' $\omega_\sigma$,
we have that
$ \omega_\sigma \big( {\cal G}_m \big)= {\GG}_{m /\QQ}$, which 
is the Mumford-Tate group of ${\rm T}_\sigma ( \ZZ(1))$.
\par (3) If $k$ is algebraically closed, the motivic Galois group of motives 
of CM-type over $k$ is the Serre group (cf. [10] 4.8).
\vskip 0.5 true cm

\pn 2.7. {\bf Lemma-Definition} 
\pn {\it The motivic Galois group $\G( {\cal T}_{0}(k))$ of ${\cal T}_{0}(k)$
is the affine group $\QQ\,$-scheme $\gal$ endowed 
with its action on itself by inner automorphisms. We denote it by $\galm$
In particular, for any fibre functor $\omega$ over $\spec(\QQ)$ of
 ${\cal T}_{0}(k),$ the affine group scheme 
$\omega(\galm)= {\underline {\rm Aut}}^{\otimes}_{\spec(\QQ)}(\omega)$ 
is canonically isomorphic to $\gal$ }
\vskip 0.5 true cm

\pn PROOF: Since $ {\cal T}_{0}(k) \cong {\rm Rep}_{\QQ}(\gal),$ this lemma
is an immediate consequence of remark 1.4 (2). 
\vskip 0.5 true cm

\pn 2.8. {\bf Lemma}
 \pn {\it
\par (i) The Tannakian subcategory ${\cal T}_{ 0}(k)$ of $MR(k)$ is 
equivalent  (as tensor category) to the Tannakian subcategory
$ {\Gr}^W_0 {\cal T}_{ 1}(k)$.
\par (ii) We have the following anti-equivalence of tensor categories }

$$\eqalign{ {\cal T}_{ 0}(k) \otimes \langle \ZZ(1) \rangle^\otimes
 & \longrightarrow {\Gr}^W_{-2} {\cal T}_{ 1}(k) \cr 
X \otimes \ZZ(1) & \longmapsto X^{\du}(1).\cr }$$

\pn PROOF:
\pn (i) It is a consequence of (2.2.1) and of the fact that, as observed in 2.5,  in the Tannakian category $MR(k)$
 we lose the integral structures.

\pn (ii) According to (i), 
we can view an object $X$ of ${\cal T}_{ 0}(k)$ as the character group of a torus $T$  defined over $k$. 
The dual $X^{\du}$ of $X$ in the Tannakian category $MR(k),$ can be identified with the cocharacter group of $T$ which can be written, according to our notation, 
as $X^{\du}(1).$
The anti-equivalence between the category of character groups and the 
category of cocharacter groups furnishes the desired anti-equivalence.
\vskip 0.5 true cm

\pn 2.9. {\bf Corollary} 
{\it 
\par (i) $\G( {\Gr}^W_0{\cal T}_{1}(k))= \galm$,
\par (ii) $\G( {\Gr}^W_{-2}{\cal T}_{1}(k))= i_1 \galm \times i_2 {\cal G}_{m}$,
\pn  where 
  $ i_1: {\cal T}_{0}(k)={\cal T}_{0}(k) \otimes {\rm Vec} (\QQ)
\longrightarrow {\cal T}_{ 0}(k) \otimes \langle \ZZ(1) \rangle^\otimes$
 and  $ i_2: \langle \ZZ(1) \rangle^\otimes = {\rm Vec} (\QQ) \otimes
 \langle \ZZ(1) \rangle^\otimes \longrightarrow
 {\cal T}_{ 0}(k) \otimes \langle \ZZ(1) \rangle^\otimes$ 
identify respectively 
${\cal T}_{ 0}(k)$ and $\langle \ZZ(1) \rangle^\otimes$ with full subcategories of ${\cal T}_{ 0}(k) \otimes \langle \ZZ(1)\rangle^\otimes$ . (In the following, we will avoid the symbols $i_1$ and $i_2$).}
\vskip 0.5 true cm

\pn PROOF: 
\pn (i) Consequence of (2.8.1).
\pn (ii) From 2.8 (ii) and (2.40.5) [10], we have that
$$\eqalign{\G \big( {\Gr}^W_{-2}{\cal T}_{1}(k)\big)
=&\G\big({\cal T}_{ 0}(k) \otimes \langle \ZZ(1)\rangle^\otimes \big) \cr
=& i_1 \G\big({\cal T}_{ 0}(k)\big) \times i_2 {\cal G}_{m}=
i_1 \galm \times i_2 {\cal G}_{m}.\cr}$$
\vskip 0.5 true cm

\pn 2.10. REMARKS: 
\par (1) The motivic Galois group $\G( {\cal T}_{0}(\ok))$ is the affine
 group $ {\cal T}_{0}(\ok)$-scheme 
\pn $\sp (1_{{\cal T}_{0}(\ok)})$ defined by the 
affine group  $\QQ\,$-scheme $\spec(\QQ)$ (cf. (2.2.2) and 1.4 (2)).
\par (2) Since the category ${\Gr}^W_{-2}{\cal T}_{1}(\ok)$ is equivalent 
to the Tannakian category generated by the torus $\ZZ(1)$, the motivic Galois group 
$\G({\Gr}^W_{-2}{\cal T}_{1}(\ok))$ is  ${\cal G}_m $.
\par (3) In the category of affine group  
$ {\cal T}_{0}(k)$-schemes, there are two $ {\cal T}_{0}(k)$-schemes
defined by the Galois group $\gal$:
\par - the affine group $ {\cal T}_{0}(k)$-scheme $\galm$ which is the affine 
group $\QQ\,$-scheme $\gal$ endowed with its action on itself
 by inner automorphisms. It is the fundamental group of the Tannakian category 
$ {\cal T}_{0}(k)$ of Artin motives.
\par - the affine group $ {\cal T}_{0}(k)$-scheme $\gal$  which is the affine 
group $\QQ\,$-scheme $\gal$ endowed with the trivial action of $\gal$
(cf. 1.2 (2)). It is a $\QQ\,$-scheme
viewed as a $ {\cal T}_{0}(k)$-scheme.
\pn Same remark for the affine group $\langle \ZZ(1) \rangle^\otimes$-schemes ${\cal G}_{m}$ and ${\GG}_m.$
\vskip 0.5 true cm

\pn 2.11. The weight filtration $W_*$ on objects of ${\cal T}_{1}(k)$ induces an increasing filtration, always denoted by $W_*,$ on the motivic Galois group ${\G}({\cal T}_{ 1}(k))$ of ${\cal T}_{ 1}(k)$ (cf. [11] Chapitre IV \S 2). We describe this filtration through the action (1.3.3) of ${\G}({\cal T}_{ 1}(k))$ on the generators of ${\cal T}_{1}(k).$ For each 1-motive $M$ over $k$,

\vskip 0.5 true cm
\pn $ \quad  W_{0}({\G}({\cal T}_{1}(k)))={\G}({\cal T}_{1}(k)) $
\vskip 0.3 true cm

\pn $ \quad  W_{-1}({\G}({\cal T}_{1}(k)))=\{ g \in {\G}({\cal T}_{1}(k)) \, \, \vert \, \,
(g - id)M \subseteq   W_{-1}(M),$
\pn $ (g - id) W_{-1}(M) \subseteq W_{-2}(M), (g - id) W_{-2}(M)=0 \} , $
\vskip 0.3 true cm

\pn $  \quad  W_{-2}({\G}({\cal T}_{1}(k)))=\{ g \in {\G}({\cal T}_{1}(k)) \, \, \vert \, \,
(g - id) M \subseteq W_{-2}(M), $
\pn $ (g - id)  W_{-1}(M) =0\}, $
\vskip 0.3 true cm

\pn $ \quad  W_{-3}({\G}({\cal T}_{1}(k)))=0.$
\vskip 0.5 true cm

\pn The step $ W_{-1}({\G}({\cal T}_{1}(k)))$ is an affine unipotent group sub-${\cal T}_{1}(k)$-scheme
of the motivic Galois group ${\G}({\cal T}_{ 1}(k))$.
According to 2.9 and to the motivic analogue of [3] \S 2.2, 
${\rm Gr}_{0}^{W}({\G}({\cal T}_{1}(k)))$ acts through $\galm$ on ${\rm Gr}_{0}^{W}{\cal T}_{1}(k)$,
and through $ \galm \times  {\cal G}_m$
on ${\rm Gr}_{-2}^{W}{\cal T}_{ 1}(k)$. Moreover, its image in the automorphisms group of
${\rm Gr}_{-1}^{W}{\cal T}_{1}(k)$ is the motivic Galois group of the abelian part
 ${\rm Gr}_{-1}^{W}{\cal T}_{1}(k).$ Therefore, ${\rm Gr}_{0}^{W}({\G}({\cal T}_{1}(k)))$
is an affine reductive group sub-${\cal T}_{1}(k)$-scheme
of ${\G}({\cal T}_{ 1}(k))$ and $W_{-1}({\G}({\cal T}_{1}(k)))$ is the unipotent radical of ${\G}({\cal T}_{1}(k))$.
\vskip 0.5 true cm

\pn 2.12. Consider the base extension functor

$$\eqalign{E: {\cal T}_{1}(k)& \longrightarrow {\cal T}_{1}(\ok)\cr
                   M & \longmapsto M \otimes_{k} {\ok}. \cr }\leqno(2.12.1)$$

\pn According to (2.2.2), the kernel of $E$ is ${\cal T}_{0}(k)$,
i.e. we have the exact sequence 

$$ 0 \longrightarrow {\cal T}_{0}(k)\longrightarrow  {\cal T}_{1}(k)
{\buildrel E \over \longrightarrow}  {\cal T}_{1}(\ok).\leqno(2.12.2)$$

\pn If $M$ is an object of $ {\cal T}_{1}(\ok)$, it can be written as a 
subquotient of $M' \otimes_{k} \ok$ for some object $M'$ of 
${\cal T}_{1}(k):$ in fact, for $M'$ we can take the restriction of scalars 
 ${\rm Res}_{k'/k}M_0$ with 
$M_0$ a model of $M$ over a finite extension $k'$ of $k$. 
By 1.5, this means that the corresponding 
morphism of affine group ${\cal T}_{1}(\ok)$-schemes

$$e: \G({\cal T}_{1}(\ok)) \hookrightarrow E \, \G( {\cal T}_{1}(k))
\leqno(2.12.3)$$

\pn is a closed immersion. 

\par Finally, denote by  

$$\eqalign{
{\rm Res}_{\ok/k}: {\cal T}_{1}(\ok)& \longrightarrow {\cal T}_{1}(k)\cr
        M & \longmapsto {\rm Res}_{\ok/k}(M)   \cr }\leqno(2.12.4)$$

\pn the functor ``restriction of scalars''. 
\vskip 0.5 true cm

\pn {\bf 3. Tannakian subcategories of ${\cal T}_{ 1}(k)$.}
\vskip 0.5 true cm

\pn 3.1. Consider the following diagram of inclusions of Tannakian categories

$$\matrix{
{\Gr}_0^W {\cal T}_{1}(k)={\cal T}_{0}(k)  & & & & \cr
  &{\buildrel I_0 \over \searrow} & & & \cr
{\Gr}_{-1}^W {\cal T}_{1}(k) &{\buildrel I_{-1} \over \rightarrow} &{\Gr}_*^W {\cal T}_{1}(k) & {\buildrel I \over \rightarrow} & {\cal T}_{1}(k) \cr
 &{\buildrel I_{-2} \over \nearrow} & & & \cr
{\Gr}_{-2}^W {\cal T}_{1}(k) & & & & \cr}$$

\pn and denote by 
${\cal I}_j:{\Gr}_{j}^W {\cal T}_{1}(k)\longrightarrow {\cal T}_{1}(k)$
the inclusion $I \circ I_j$ for $j=0,-1,-2$.
\pn According to 1.5, to this diagram corresponds the diagram of faithfully flat morphisms of affine group ${\cal T}_{1}(k)$-schemes

$$\matrix{
 & & & &{\cal I}_0 \galm  \cr
  & & &{\buildrel Ii_{o} \over \nearrow} & \cr
\G(  {\cal T}_{1}(k))&{\buildrel i \over \rightarrow} &I \G({\Gr}_*^W {\cal T}_{1}(k)) & {\buildrel Ii_{-1} \over \rightarrow} &{\cal I}_{-1}\G({\Gr}_{-1}^W {\cal T}_{1}(k))  \cr
 & & & {\buildrel Ii_{-2} \over \searrow}& \cr
 & & & &{\cal I}_{-2} (\galm \times {\cal G}_m) \cr}$$

\pn where $Ii_{j}: I \G({\Gr}_*^W {\cal T}_{1}(k)) \longrightarrow 
{\cal I}_{j}\G({\Gr}_{j}^W {\cal T}_{1}(k))$ is the 
morphism of affine group ${\cal T}_{1}(k)$-schemes defined by the 
morphism of affine group ${\Gr}_*^W{\cal T}_{1}(k)$-schemes
 $i_{j}:  \G({\Gr}_*^W {\cal T}_{1}(k)) \longrightarrow 
{ I}_{j}\G({\Gr}_{j}^W {\cal T}_{1}(k))$ corresponding to the inclusion $I_j.$
\pn For $j=0,-1,-2,$ denote by 

$$\iota_j:\G(  {\cal T}_{1}(k)) \longrightarrow {\cal I}_{j}\G({\Gr}_{j}^W {\cal T}_{1}(k)) $$

\pn the faithfully flat
morphism of affine group ${\cal T}_{1}(k)$-schemes corresponding to the 
inclusion ${\cal I}_j$.  
In particular, by 1.11 we have that $\iota_j= Ii_j \circ i$.

\par Always according to 1.5, the functor ``take the graduated'' ${\Gr}_*^W: {\cal T}_{1}(k) \longrightarrow {\Gr}_*^W {\cal T}_{1}(k)$ 
corresponds to the closed immersions of affine group 
${\Gr}_*^W{\cal T}_{1}(k)$-schemes

$$ {\gr}_{*}^W:\G( {\Gr}_*^W {\cal T}_{1}(k)) \longrightarrow {\Gr}_{*}^W \G({\cal T}_{1}(k)),$$  

\pn which identifies the motivic Galois group of ${\Gr}_*^W {\cal T}_{1}(k)$ with the quotient 
${\Gr}_0^W$ of the motivic Galois group of ${\cal T}_{1}(k).$

\pn Consider now the following commutative diagram of inclusion of Tannakian categories

$$\matrix{
 {\Gr}_0^W {\cal T}_{1}(k)={\cal T}_{0}(k)  & {\buildrel J_3 \over \rightarrow}& W_0/W_{-2} {\cal T}_{1}(k) & & \cr
  &{\buildrel J_4 \over \nearrow} & &{\buildrel J_1 \over \searrow}  & \cr
{\Gr}_{-1}^W {\cal T}_{1}(k) & & {\buildrel {\cal I}_{-1} \over \rightarrow}&& {\cal T}_{1}(k) \cr
 &{\buildrel J_5 \over \searrow} & & {\buildrel J_2 \over \nearrow} & \cr
 {\Gr}_{-2}^W {\cal T}_{1}(k)&{\buildrel J_6\over \rightarrow} & W_{-1} {\cal T}_{1}(k)& & \cr
}$$

\pn where $J_1 \circ J_3 = {\cal I}_0$ and $J_2 \circ J_6= {\cal I}_{-2}$.  
By 1.5, to this diagram corresponds the commutative diagram of faithfully flat morphisms of affine group
${\cal T}_{1}(k)$-schemes 

$$\matrix{
 & & J_1\G(W_0/W_{-2} {\cal T}_{1}(k)) & {\buildrel J_1j_3 \over \rightarrow} &{\cal I}_0 \galm  \cr
  & {\buildrel j_1 \over \nearrow}& &{\buildrel J_1j_4 \over \searrow} & \cr
\G(  {\cal T}_{1}(k))& &{\buildrel \iota_{-1} \over \rightarrow} & &{\cal I}_{-1}\G({\Gr}_{-1}^W {\cal T}_{1}(k))  \cr
 & {\buildrel j_2 \over \searrow}& & {\buildrel J_2j_5 \over \nearrow}& \cr
 & &J_2 \G( W_{-1} {\cal T}_{1}(k)) & {\buildrel J_2j_6 \over \rightarrow} &{\cal I}_{-2}
 (\galm \times {\cal G}_m) \cr}\leqno(3.1.1)$$

\pn where by 1.11 we have that $\iota_0= J_1j_3 \circ j_1$ and $ \iota_{-2}=J_2j_6 \circ j_2.$ 
\vskip 0.5 true cm

\pn 3.2. According to 1.7, the Tannakian subcategories $ {\cal T}_{0}(k), {\rm  Gr}_{-1}^W {\cal T}_{1}(k),$
\pn $ {\rm  Gr}_{-2}^W {\cal T}_{1}(k), W_0/W_{-2} {\cal T}_{1}(k), W_{-1} {\cal T}_{1}(k)$ are characterized by the following af\-fi\-ne 
\pn group 
sub-${\cal T}_{1}(k)$-schemes of $\G({\cal T}_{1}(k))$.

$$\eqalign{
 H_{{\cal T}_{1}(k)}({\rm  Gr}_{0}^W {\cal T}_{1}(k)) &=\ker \big[ \G({\cal T}_{1}(k))\,  {\buildrel \iota_0 \over \longrightarrow} \,{\cal I}_0 \galm \big],\cr
H_{{\cal T}_{1}(k)}({\rm  Gr}_{-1}^W {\cal T}_{1}(k)) &=\ker \big[ \G({\cal T}_{1}(k))\,  {\buildrel \iota_{-1} \over \longrightarrow} \,{\cal I}_{-1} \G({\Gr}_{-1}^W {\cal T}_{1}(k)) \big],\cr
H_{{\cal T}_{1}(k)}({\rm  Gr}_{-2}^W {\cal T}_{1}(k)) &=\ker \big[ \G({\cal T}_{1}(k))\,  {\buildrel \iota_{-2} \over \longrightarrow} \,{\cal I}_{-2} ( \galm \times {\cal G}_m) \big],\cr
H_{{\cal T}_{1}(k)}(W_0/W_{-2} {\cal T}_{1}(k))& =\ker \big[ \G({\cal T}_{1}(k))\, {\buildrel j_1 \over \longrightarrow} \, J_1 \G(W_0/W_{-2} {\cal T}_{1}(k))\big],\cr
H_{{\cal T}_{1}(k)}(W_{-1} {\cal T}_{1}(k)) & =\ker \big[ \G({\cal T}_{1}(k)) \, {\buildrel j_2 \over \longrightarrow} \, J_2 \G(W_{-1} {\cal T}_{1}(k)) \big].\cr}$$

\pn These affine group sub-${\cal T}_{1}(k)$-schemes are the motivic generalization of the algebraic $\QQ \, $-groups introduced in [1] \S 2. 
Because of the commutativity of the diagram (3.1.1), we have the inclusions 

$$ \eqalign{
H_{{\cal T}_{1}(k)}(W_{-1} {\cal T}_{1}(k))&\subseteq  
H_{{\cal T}_{1}(k)}({\rm  Gr}_{-1}^W {\cal T}_{1}(k)) \cap H_{{\cal T}_{1}(k)}({\rm  Gr}_{-2}^W {\cal T}_{1}(k)),\cr
H_{{\cal T}_{1}(k)}(W_0/W_{-2} {\cal T}_{1}(k))&\subseteq 
H_{{\cal T}_{1}(k)}({\rm  Gr}_{0}^W {\cal T}_{1}(k))\cap H_{{\cal T}_{1}(k)}({\rm  Gr}_{-1}^W {\cal T}_{1}(k)).\cr}$$

\pn By lemma 2.8, the Tannakian category ${\rm  Gr}_{0}^W {\cal T}_{1}(k)$ of Artin motives is a Tannakian subcategory of ${\rm  Gr}_{-2}^W {\cal T}_{1}(k)$ and therefore according to 1.8 we have that
$$ H_{{\cal T}_{1}(k)}({\rm  Gr}_{0}^W {\cal T}_{1}(k))\supseteq  H_{{\cal T}_{1}(k)}({\rm  Gr}_{-2}^W {\cal T}_{1}(k)).\leqno(3.2.1)$$

\pn Moreover the Cartier duality furnishes the anti-equivalence of tensor categories
$$\eqalign{
W_0/W_{-2} {\cal T}_{1}(k) \, & \longrightarrow \, W_{-1} {\cal T}_{1}(k)\cr 
M & \longmapsto M^*.\cr}$$

\pn We denote both these categories by ${\widetilde W}{\cal T}_{1}(k)$. In particular we have that
$$H_{{\cal T}_{1}(k)}(W_0/W_{-2} {\cal T}_{1}(k))=H_{{\cal T}_{1}(k)}(W_{-1} {\cal T}_{1}(k))=H_{{\cal T}_{1}(k)} ({\widetilde W}{\cal T}_{1}(k)).\leqno(3.2.2)$$

\pn With these notation, we can state the motivic generalization of lemma [1] 2.2
to all motives of niveau $\leq 1$.
\vskip 0.5 true cm

\pn 3.3. {\bf Lemma}
\pn (i) $ W_{-1}(\G({\cal T}_{1}(k)))= \cap_{i=-1,-2}  H_{{\cal T}_{1}(k)}({\rm  Gr}_{i}^W {\cal T}_{1}(k)) , $  
\pn (ii) $ W_{-2}(\G({\cal T}_{1}(k)))= H_{{\cal T}_{1}(k)}({\widetilde W}{\cal T}_{1}(k)). $
\vskip 0.5 true cm

\pn PROOF: Let $M$ be a 1-motive over $k$.
 We remark  that

$$\eqalign{ g \in H_{{\cal T}_{1}(k)}({\rm  Gr}_{0}^W {\cal T}_{1}(k)) 
& \Longleftrightarrow g_{\vert {\rm Gr}_{0}^{W} (M) } =id 
 \Longleftrightarrow  (g-id ) M \subseteq W_{-1}(M) \cr
          g \in H_{{\cal T}_{1}(k)}({\rm  Gr}_{-1}^W {\cal T}_{1}(k))
& \Longleftrightarrow  g_{\vert {\rm Gr}_{-1}^{W}( M)} =id 
 \Longleftrightarrow (g-id ) W_{-1}(M) \subseteq W_{-2}(M) \cr
             g \in H_{{\cal T}_{1}(k)}({\rm  Gr}_{-2}^W {\cal T}_{1}(k)) 
& \Longleftrightarrow  g_{\vert {\rm Gr}_{-2}^{W} (M)} =id 
 \Longleftrightarrow (g-id ) W_{-2}(M) =0 \cr
           g \in H_{{\cal T}_{1}(k)}(W_0/W_{-2} {\cal T}_{1}(k)) 
& \Longleftrightarrow  g_{\vert W_0/ W_{-2}(M) } =id 
 \Longleftrightarrow (g-id ) M  \subseteq  W_{-2}(M) \cr
             g \in H_{{\cal T}_{1}(k)}(W_{-1} {\cal T}_{1}(k)) 
& \Longleftrightarrow  g_{\vert  W_{-1}(M) } =id  \Longleftrightarrow
 (g-id )   W_{-1}(M)  =0 . \cr}$$
 
\pn The result now follows from (3.2.1), (3.2.2) and from the definitions of 
\pn $W_{-1}(\G({\cal T}_{1}(k)))$ and of $W_{-2}(\G({\cal T}_{1}(k)))$.
\vskip 0.5 true cm

\pn 3.4. {\bf Theorem} 
\pn {\it We have the following diagram of morphisms of 
affine group sub-${\cal T}_{1}(k)$-schemes}

$$\matrix{
{\scriptstyle 0} &\rightarrow &{\scriptstyle {\rm Res}_{\ok/k} \G({\cal T}_{1}(\ok))} &\rightarrow &{\scriptstyle \G({\cal T}_{1}(k))} & {\buildrel \iota_0 \over  \rightarrow } &{\scriptstyle {\cal I}_0 \galm} &\rightarrow &{\scriptstyle 0} \cr
&&\uparrow &&||&&\uparrow& \cr
{\scriptstyle 0} &\rightarrow &{\scriptstyle {\rm Res}_{\ok/k} H_{{\cal T}_{1}(\ok)}(\langle \ZZ(1) \rangle^\otimes) } &\rightarrow &{\scriptstyle \G({\cal T}_{1}(k))} & {\buildrel \iota_{-2} \over  \rightarrow } &{\scriptstyle {\cal I}_{-2} (\galm \times {\cal G}_m) } &\rightarrow &{\scriptstyle 0} \cr
&&\uparrow &&||&&\uparrow& \cr
{\scriptstyle 0} &\rightarrow &{\scriptstyle W_{-1} \G({\cal T}_{1}(k))} &\rightarrow &{\scriptstyle \G({\cal T}_{1}(k))} & {\buildrel i \over  \rightarrow } &{\scriptstyle I\G({\Gr}_*^W {\cal T}_{1}(k))} &\rightarrow &{\scriptstyle 0} \cr
&&\uparrow &&||&&\uparrow& \cr
{\scriptstyle 0} &\rightarrow &{\scriptstyle W_{-2} \G({\cal T}_{1}(k))} &\rightarrow &{\scriptstyle \G({\cal T}_{1}(k))} & {\buildrel j_1\over  \rightarrow } &{\scriptstyle J_1\G({\widetilde W}{\cal T}_{1}(k))} &\rightarrow &{\scriptstyle 0} \cr
}$$

\pn {\it where all horizontal short sequences are exact and where the vertical arrows on the left are inclusions and those on the right are surjections.}
\vskip 0.5 true cm

\pn PROOF: We will prove the exactness of the four horizontal short sequences applying theorem 1.7 to the following 
Tannakian subcategories of ${\cal T}_{1}(k)$: ${\cal T}_{0}(k),
{\widetilde W}{\cal T}_{1}(k),$
\pn ${\Gr}_{-2}^W {\cal T}_{1}(k)$ and ${\Gr}_*^W {\cal T}_{1}(k)$.

\pn By [8] 2.20 (e) the composite functor ${\rm Res}_{\ok/k} \circ E$ is the multiplication by $[\ok :k].$ Therefore, since we work modulo isogenies, according (2.12.3) the affine group ${\cal T}_{1}(k)$-scheme 
${\rm Res}_{\ok/k} \G( {\cal T}_{1}(\ok))$ is a 
sub-${\cal T}_{1}(k)$-scheme of $\G( {\cal T}_{1}(k)).$ 
Now consider the inclusion ${\cal I}_0: 
 {\Gr}_0^W {\cal T}_{1}(k) \longrightarrow  {\cal T}_{1}(k)$ (cf. 3.1). 
Since Artin motives are the kernel of the base extension functor (cf. (2.12.2)),
the objects of ${\cal T}_{0}(k)$ are 
exactly those objects of $ {\cal T}_{1}(k)$ on which the
sub-${\cal T}_{1}(k)$-scheme ${\rm Res}_{\ok/k} \G({\cal T}_{1}(\ok))$ 
of $\G( {\cal T}_{1}(k))$ acts trivially. Hence by 1.7, we have the exact sequence of affine group 
${\cal T}_{1}(k)$-schemes

$$0 \longrightarrow  {\rm Res}_{\ok/k} \G({\cal T}_{1}(\ok)) \longrightarrow
\G({\cal T}_{1}(k)) {\buildrel \iota_0 \over \longrightarrow} 
{\cal I}_0 \galm \longrightarrow 0.$$

\pn As in 3.1 denote by ${\cal I}_{-2}$ the inclusion of the Tannakian subcategory ${\Gr}_{-2}^W {\cal T}_{1}(k)$ in the Tannakian category ${\cal T}_{1}(k)$.
 As observed in 2.10 (2), $ {\Gr}_{-2}^W {\cal T}_{1}(\ok)$ is equivalent as tensor category to the Tannakian subcategory $\langle \ZZ(1) \rangle^\otimes$ of ${\cal T}_{1}(\ok)$ generated by the $k$-torus $\ZZ(1)$. Hence the objects of $ {\Gr}_{-2}^W {\cal T}_{1}(k)$ are exactly those objects
 of ${\cal T}_{1}(k)$ on which, after extension of scalars, the sub-${\cal T}_{1}(\ok)$-scheme
$ H_{{\cal T}_{1}(\ok)}(\langle \ZZ(1) \rangle^\otimes)$ of $\G({\cal T}_{1}(\ok))$ acts trivially.
Applying 1.7, we have the exact sequence of affine group ${\cal T}_{1}(k)$-schemes
 
 $$0 \longrightarrow  {\rm Res}_{\ok/k} H_{{\cal T}_{1}(\ok)}(\langle \ZZ(1) \rangle^\otimes) \longrightarrow
\G({\cal T}_{1}(k)) {\buildrel \iota_{-2} \over \longrightarrow} 
 {\cal I}_{-2} ( \galm \times {\cal G}_m ) \longrightarrow 0$$
 
\pn Consider now the inclusion $ i: {\Gr}_*^W {\cal T}_{1}(k)
 \longrightarrow  {\cal T}_{1}(k)$ (cf. 3.1). 
 Since the objects of ${\Gr}_*^W {\cal T}_{1}(k)$ are 
exactly those objects of $ {\cal T}_{1}(k)$ on which the unipotent 
radical $W_{-1} \G({\cal T}_{1}(k)) $ of $\G({\cal T}_{1}(k)) $ acts
trivially, the affine group ${\cal T}_{1}(k)$-scheme characterizing 
${\Gr}_*^W {\cal T}_{1}(k)$ is

$$H_{{\cal T}_{1}(k)} ({\Gr}_*^W {\cal T}_{1}(k))= W_{-1} \G({\cal T}_{1}(k)) $$

\pn (cf. [6] 6.7 (ii)). In particular, we have the exact sequence of affine group 
${\cal T}_{1}(k)$-schemes

$$0 \longrightarrow  W_{-1} \G({\cal T}_{1}(k)) \longrightarrow
\G({\cal T}_{1}(k)) {\buildrel i \over \longrightarrow} 
I\G({\Gr}_*^W {\cal T}_{1}(k)) \longrightarrow 0$$

\pn and as expected the motivic Galois group $\G({\Gr}_*^W {\cal T}_{1}(k))$ is isomorph to
\pn $\G({\cal T}_{1}(k))/ W_{-1} \G({\cal T}_{1}(k)).$

\pn As in 3.1 let $J_1: 
 {\widetilde W}{\cal T}_{1}(k) \longrightarrow  {\cal T}_{1}(k)$ be the inclusion 
 of the Tannakian subcategory ${\widetilde W}{\cal T}_{1}(k)$ in the 
 Tannakian category ${\cal T}_{1}(k).$
 According to 3.3 (ii) and 1.7 we have the exact sequence of affine group ${\cal T}_{1}(k)$-schemes
 
 $$0 \longrightarrow  W_{-2} \G({\cal T}_{1}(k)) \longrightarrow
\G({\cal T}_{1}(k)) {\buildrel j_1 \over \longrightarrow} 
J_1 \G({\widetilde W} {\cal T}_{1}(k)) \longrightarrow 0.$$

\pn Finally, in order to prove that the left vertical arrows are inclusions and that 
the right vertical arrows are surjections, it is enough to apply the lemma 1.8.
\vskip 0.5 true cm

\pn 3.5. As corollary, we get the motivic version of [5] II 6.23 (a), (c)
 and [8] 4.7 (c), (e). But before to state this corollary, we recall some facts:
\par If $F_1,F_2:{\cal T}_{1}(\ok) \longrightarrow {\cal T}_{1}(\ok)$ 
are two functors, we define $ {\underline {\rm Hom}}^{\otimes}(F_1,F_2)$ to be 
the functor which associates to each  ${\cal T}_1(\ok)$-scheme
${\rm Sp}(B)$, the set of morphisms of
 $\otimes$-functors from $(F_1)_{{\rm Sp}(B)}:X \longmapsto F_1(X)\otimes B$
to $(F_2)_{{\rm Sp}(B)}:X \longmapsto F_2(X)\otimes B$ ($(F_1)_{{\rm Sp}(B)}$ and $(F_2)_{{\rm Sp}(B)}$
are $\otimes$-functors from ${\cal T}_{1}(\ok)$ to the category of modules over ${\rm Sp}(B)$).
\par Moreover, each element $\tau$ of $\gal$ defines a functor

$$\tau:{\cal T}_{1}(\ok) \longrightarrow {\cal T}_{1}(\ok)$$

\pn in the following way: since as observed in 2.12, the category ${\cal T}_{1}(\ok)$ is generated 
by motives of the form $E(M)$ with $M \in {\cal T}_{1}(k),$ it is enough to define $\tau E(M).$
We put $\tau E(M)= M \otimes_k \tau \ok.$
\vskip 0.5 true cm

\pn 3.6. {\bf Corollary} 
\pn {\it (i) We have the following diagram of morphisms of 
affine group sub-${\cal T}_{1}(k)$-schemes in which all the short sequences 
are exact:}
$$\matrix{
 & & & &{\scriptstyle 0} & & & & \cr
 & & & &\downarrow & & & & \cr
 & & & &{\scriptstyle W_{-1} \G({\cal T}_{1}(k))} &&& & \cr
 & & & &\downarrow & & & & \cr
{\scriptstyle 0} &\rightarrow &{\scriptstyle {\rm Res}_{\ok/k} \G({\cal T}_{1}(\ok))} &\rightarrow &{\scriptstyle \G({\cal T}_{1}(k))} & {\buildrel \iota_0 \over  \rightarrow } &{\scriptstyle {\cal I}_0 \galm} &\rightarrow &{\scriptstyle 0} \cr
 & &\downarrow & & {\scriptstyle i} \downarrow & &\Vert & & \cr
{\scriptstyle 0} &\rightarrow &{\scriptstyle I {\rm Res}_{\ok/k} \G({\Gr}_*^W {\cal T}_{1}(\ok))} &\rightarrow &{\scriptstyle I\G({\Gr}_*^W {\cal T}_{1}(k))} & {\buildrel I i_0 \over  \rightarrow } &{\scriptstyle I I_0 \galm} &\rightarrow &{\scriptstyle 0} \cr
 &  & \downarrow& &\downarrow & & & & \cr
 & & {\scriptstyle 0}& &{\scriptstyle 0} & & & & \cr
}$$

\pn {\it (ii) The morphism $ I{\gr}_{*}^W:I \G({\Gr}_*^W {\cal T}_{1}(k)) \longrightarrow  \G({\cal T}_{1}(k))$ of affine group ${\cal T}_{1}(k)$-schemes 
is a section of $i.$
\pn (iii) For any $\tau \in \galm,$ $\iota_0^{-1}(\tau)= {\underline {\rm Hom}}^{\otimes}({\rm Id},\tau \circ {\rm Id} ),$ regarding $ {\rm Id}$ and $ \tau \circ {\rm Id}$ as functors on ${\cal T}_{1}(\ok)$. In an analogous way, $I i_0^{-1}(\tau)= {\underline {\rm Hom}}^{\otimes}({\rm Id},\tau \circ {\rm Id}),$ regarding $ {\rm Id}$ and $ {\rm Id}\circ \tau$ as functors on 
$ {\Gr}_*^W {\cal T}_{1}(\ok)$.}
\vskip 0.5 true cm

\pn PROOF: $(i)$ We have only to prove the exactness of the second horizontal short sequence. 
In order to do this, we apply theorem 1.7 to the Tannakian category 
$  {\Gr}_0^W {\cal T}_{1}(k)$  viewed this time as
 subcategory of $  {\Gr}_*^W {\cal T}_{1}(k)$. 
Consider the inclusion $I_0:  {\Gr}_0^W {\cal T}_{1}(k) \longrightarrow 
 {\Gr}_*^W {\cal T}_{1}(k)$ (cf 3.1). By (2.12.2), Artin  motives are exactly the kernel of the base
 extension functor and therefore 
 
 $$H_{{\Gr}_*^W {\cal T}_{1}(k)} ({\Gr}_0^W {\cal T}_{1}(k))= {\rm Res}_{\ok/k} \G({\Gr}_*^W {\cal T}_{1}(\ok)).$$

\pn $(ii)$ Since

$$I \circ {\Gr}_*^W = id: {\cal T}_{1}(k)\, {\buildrel {\Gr}_*^W \over \longrightarrow }
 {\Gr}_*^W {\cal T}_{1}(k) \, {\buildrel I \over \longrightarrow }\, {\cal T}_{1}(k),$$

\pn from 1.11 we have that

$$ I{\gr}_*^W \circ i= id: \G( {\cal T}_{1}(k)) \, {\buildrel i \over \longrightarrow}\,
I \G( {\Gr}_*^W {\cal T}_{1}(k))\, {\buildrel I{\gr}_*^W  \over \longrightarrow}
\, \G( {\cal T}_{1}(k)).$$

\pn $(iii)$  By [7] 8.11, the fundamental group
$\G ({\cal T}_1(k))$ represents the functor 
\pn $ {\underline {\rm Aut}}^{\otimes}({\rm Id})$ 
which associates to each  ${\cal T}_1(k)$-scheme
${\rm Sp}(B)$ the group of automorphisms of
 $\otimes$-functors of the functor

$$\eqalign{ {\rm Id}_{{\rm Sp}(B)}: {\cal T}_1(k)& \longrightarrow  \{
 {\rm modules~ over~} {\rm Sp}(B) \} \cr 
X& \longmapsto X \otimes B.\cr }$$

\pn Hence if $g$ is an element of 
$\G( {\cal T}_{1}(k))({\sp} B)= {\underline {\rm Aut}}^{\otimes}({\rm Id})({\sp} B),$ for each pair of objects $M$ and $N$ of ${\cal T}_{1}(k)$ and 
for each morphism $f:M \longrightarrow N$ of ${\cal T}_{1}(k)$, we have the 
commutative diagram
 
$$\matrix{ M \otimes B & {\buildrel g_M \over \longrightarrow}& M \otimes B \cr
f \otimes id_B \downarrow~~~~~~~~~~& &~~~~~~~~~~~~ \downarrow f \otimes id_B \cr
N \otimes B & {\buildrel g_N \over \longrightarrow}& N \otimes B. \cr }$$

\pn Let $M$ and $N$ be two objects of ${\cal T}_{1}(k)$. 
Since ${\Hom}_{{\cal T}_{1}(\ok)} (E(M),E(N))$
is an object of ${\rm Rep}_{\QQ}(\gal)$, it can be regarded as an Artin motive over $k$. Moreover, 
 the elements of ${\Hom}_{{\cal T}_{1}(k)} (M,N)$ are exactly the elements of 
 ${\Hom}_{{\cal T}_{1}(\ok)} (E(M),E(N))$ which are invariant under the action of $\gal$, i.e.

$${\Hom}_{{\cal T}_{1}(k)} (M,N)=\big({\Hom}_{{\cal T}_{1}(\ok)} (E(M),E(N))\big)^{\gal}.$$

\pn Let $g$ be an element of $\G( {\cal T}_{1}(k))({\sp} B)= {\underline {\rm Aut}}^{\otimes}({\rm Id})({\sp} B),$ and let $\iota_0(g)=\tau \in \gal.$ This means that $g$ acts via $\tau$ on 
${\Hom}_{{\cal T}_{1}(\ok)} (E(M),E(N))$.
Then for any  morphism $h:E(M) \longrightarrow E(N)$ of ${\cal T}_{1}(\ok)$, we have the commutative diagram 

$$\matrix{ E(M) \otimes B & {\buildrel E(g_M) \over \longrightarrow}& E(M) \otimes B~~~~~~~~\cr
 h \otimes id_B \downarrow ~~~~~~& &~~~~~~~ \downarrow \tau h\otimes id_B \cr
E(N) \otimes B & {\buildrel E(g_N) \over \longrightarrow}&E(N) \otimes B.~~~~~~~~ \cr }\leqno(3.6.1)$$

\pn Since $M$ and $N$ are defined over $k$,
$E(M)$ and $E(N)$ are respectively isomorph to $ \tau E(M)$ and $ \tau E(N)$ and therefore
 the upper line of (3.6.1) defines a morphism 

$$  E(M) \otimes B \longrightarrow \tau E(M) \otimes B$$ 

\pn which is functorial in $E(M)$ and $B,$ and which is compatible with tensor products. Moreover 
we have already observed in 2.12 that the Tannakian category ${\cal T}_{1}(\ok)$ is generated by motives of the form $E(M)$ with $M \in {\cal T}_{1}(k)$. We can then conclude that $g$ defines an element of ${\underline {\rm Hom}}^{\otimes}({\rm Id},\tau \circ {\rm Id} ),$ regarding $ {\rm Id}$ and $ \tau \circ {\rm Id}$ as functors on ${\cal T}_{1}(\ok)$. 
\vskip 0.5 true cm

\pn {\bf 4. Specialization to a 1-motive.}
\vskip 0.5 true cm

\pn 4.1. In this section we need a more symmetric description of 1-motives: 
consider the 7-uplet $(X,Y^{\du},A,A^*, v ,v^*,\psi)$ where
\par - $X$ and $Y^{\du}$ are two group $k$-schemes, which are locally 
for the \'etale topology, constant group schemes defined by a
 finitely generated free $\ZZ$-module;
\par - $A$ and $A^*$ are two abelian varieties defined over $k$, dual to each other;
\par - $v:X \longrightarrow A$ and $v^*:Y^{\du} \longrightarrow A^*$ are two morphisms  of group schemes over $k$; and 
\par - $\psi$ is a trivialization of the pull-back $(v,v^*)^*{\cal P}_A$ by $(v,v^*)$
of the Poincar\' e biextension ${\cal P}_A$ of $(A,A^*)$.
\pn  By [4] (10.2.14), to have the data $(X,Y^{\du},A,A^*, v ,v^*,\psi)$ is equivalent to have the 1-motive $M=[X {\buildrel u \over \longrightarrow }G]$,
where $G$ is the semi-abelian variety defined by the homomorphism $v^*$ and $u$ is the lifting of $v$ determined by the trivialization $\psi.$
\vskip 0.5 true cm

\pn 4.2. Let $M=(X,Y^{\du},A,A^*, v ,v^*,\psi)$ be a 1-motive defined over $k$. It is not true that the Tannakian category generated by $W_0/W_{-2} (M)$ is equivalent to the Tannakian category generated by $ W_{-1}(M)$. Hence the lemma 3.3. must be modified in the following way

$$\eqalign{
W_{-1}(\G(M))&= \cap_{i=-1,-2}  H_{{\cal T}_{1}(k)}({\rm  Gr}_{i}^W M) ,\cr  
W_{-2}(\G(M))&= H_{{\cal T}_{1}(k)}(W_0/W_{-2} (M)) \cap H_{{\cal T}_{1}(k)}(W_{-1}(M)).\cr}\leqno(4.2.1)$$

\pn Theorem 3.4 becomes: We have the following diagram of morphisms of 
affine group sub-$\langle M \rangle^\otimes $-schemes

$$\matrix{
{\scriptstyle 0} &\rightarrow &{\scriptstyle W_{-1} \G(M)} &\rightarrow &{\scriptstyle \G(M)} & {\buildrel i \over  \rightarrow } &{\scriptstyle I\G({\Gr}_*^W M)} &\rightarrow &{\scriptstyle 0} \cr
&&\uparrow &&||&&\uparrow& \cr
{\scriptstyle 0} &\rightarrow &{\scriptstyle W_{-2} \G(M)} &\rightarrow &{\scriptstyle \G(M)} & {\buildrel j\over  \rightarrow } &{\scriptstyle J\G( W_0/W_{-2} (M)+ W_{-1}(M))} &\rightarrow &{\scriptstyle 0} \cr
&&\uparrow &&\uparrow&&\uparrow& \cr
{\scriptstyle 0} &\rightarrow &{\scriptstyle W_{-2} \G(M)} &\rightarrow &{\scriptstyle W_{-1} \G(M)} &  \rightarrow  &{\scriptstyle W_{-1} J\G( W_0/W_{-2} (M)+ W_{-1}(M))} &\rightarrow &{\scriptstyle 0} \cr
}\leqno(4.2.2)$$

\pn where the morphisms $i$ and $j$ correspond to the inclusions 
$ I:\langle {\Gr}_*^W M \rangle^\otimes \longrightarrow
 \langle M \rangle^\otimes $ and 
$J: \langle W_0/W_{-2} (M)+ W_{-1}(M) \rangle^\otimes \longrightarrow
 \langle M \rangle^\otimes,$
where all horizontal short sequences are exact and 
where all vertical arrows are inclusions, except for 
\pn $ J\G( W_0/W_{-2} (M)+ W_{-1}(M)) \longrightarrow I\G({\Gr}_*^W M)$ which is surjective.
\vskip 0.5 true cm

\pn 4.3. In [2], it is proved that if $M=(X,Y^{\du},A,A^*, v ,v^*,\psi)$ is a 1-motive over $k$, the unipotent radical of the Lie algebra 
of $\G(M)$ is the semi-abelian variety defined by the adjoint action of the 
Lie algebra $({\Gr}^W_*(W_{-1}\L \G(M)),[\,,\,])$ on itself. 
The abelian variety $B$ and the torus $Z(1)$ underlying this semi-abelian variety can be computed explicitly. We recall here briefly their construction: 

\pn The motive $E=W_{-1}( {\underline {\End}}({\Gr}_*^W M))$
is a split 1-motive whose non trivial components are the abelian variety 
$ A \otimes X^{\du}+A^* \otimes Y$ and the torus $ X^{\du} \otimes Y(1)$. 
Moreover $E$ is endowed of a Lie crochet
$[\,,\,]$ whose non trivial component 

$$( A \otimes X^{\du}+A^* \otimes Y)  \otimes (A \otimes X^{\du}+A^* \otimes Y) \longrightarrow X^{\du} \otimes Y(1).$$

\pn is defined through the morphism $A \otimes A^* \longrightarrow \ZZ(1)$ associated by [2] (1.3.1) to the Poincar\'e biextension $\cal P$ of $(A,A^*)$ by $\ZZ(1)$.
 According to corollary [2] 2.7, this Lie crochet corresponds to a 
$\Sigma - X^{\du} \otimes Y  (1)$-torsor ${\cal B}$ living over $ A \otimes X^{\du}+A^* \otimes Y.$
As proved in [2] 3.3, the 1-motives ${\Gr}_*^W M$ and ${\Gr}_*^W M^{\du}$ are Lie $(E,[\, , \,])-$modules. In particular, $E$ acts on the components $ {\Gr}_0^W M$ and 
$ {\Gr}_0^W M^{\du}$ in the following way:

 $$\eqalign{
\alpha:& (X^{\du} \otimes A) \otimes X  \longrightarrow A \cr
\beta:&   (A^* \otimes Y) \otimes Y^{\du} \longrightarrow A^*\cr 
\gamma :& (X^{\du} \otimes Y(1)) \otimes X  \longrightarrow Y(1) \cr}\leqno(4.3.1)$$

\pn These morphisms are projections defined through the evaluation maps
$ev_{X^{\du}}: X^{\du} \otimes X \longrightarrow \ZZ(0)$ and 
 $ev_{Y^{\du}}: Y^{\du} \otimes Y \longrightarrow \ZZ(0)$ of the Tanakian category 
${\cal T}_1(k)$ (cf. [7] (2.1.2)).

\pn Thank to the morphisms  
$\delta_{ X^{\du}}: \ZZ(0) \longrightarrow X \otimes X^{\du}$ and
$ev_{X}: X \otimes X^{\du} \longrightarrow \ZZ(0)$, which characterize the notion of
 duality in a Tannakian category (cf. [7] (2.1.2)), to have the morphisms 
$v: X \longrightarrow A$ and $v^*: Y^{\du} \longrightarrow A^*$ 
 is equivalent to have a $k$-rational point  $b=(b_1,b_2)$ 
of the abelian variety $ A \otimes X^{\du}+A^* \otimes Y.$
\vskip 0.3 true cm

\par - Let $B$ be the smallest abelian sub-variety (modulo isogeny) of
 $X^{\du}  \otimes A+A^* \otimes Y$ containing the point
$b=(b_1,b_2) \in X^{\du}  \otimes A (k)\times A^* \otimes Y (k) $.
\vskip 0.3 true cm

\pn The restriction  $i^*{\cal B}$ of the $\Sigma - X^{\du} \otimes Y  (1)$-torsor ${\cal B}$ by  
the inclusion $i: B \longrightarrow X^{\du}  \otimes A \times A^* \otimes Y $ is a 
 $\Sigma-X^{\du} \otimes Y(1)$-torsor over $B$.
\vskip 0.3 true cm

\par - Denote by $Z_1$ the smallest $\gal$-module of 
$X^{\du} \otimes Y$ such that the torus $Z_1(1)$ that it defines, contains  
the image of the Lie crochet $[\, ,\,]: B \otimes B \longrightarrow X^{\du} \otimes Y(1)$.
\vskip 0.3 true cm

\pn The direct image $p_*i^* {\cal B}$ of the $\Sigma-X^{\du} \otimes Y(1)$-torsor
 $i^* {\cal B}$ by the projection
 $p:X^{\du} \otimes Y(1) \longrightarrow (X^{\du} \otimes Y/ Z_1)(1)$ is a trivial
 $\Sigma-(X^{\du} \otimes Y/ Z_1)(1)$-torsor over $B$, 
i.e. $p_*i^* {\cal B}= B \times (X^{\du} \otimes Y/ Z_1)(1)$. We denote by  
$\pi: p_*i^* {\cal B} \longrightarrow  (X^{\du} \otimes Y/ Z_1)(1)$
the canonical projection.
\pn By [2] 3.6, to have the lifting $u: X \longrightarrow G$ of the morphism 
$v:X \longrightarrow A$ is equivalent to have a point  
 $\widetilde b$ in the fibre of $ {\cal B}$ over $b$. 
We denote again by $\widetilde b$ the points of $i^* {\cal B}$
and of $p_*i^* {\cal B}$ over the point $b$ of $B$.
\vskip 0.3 true cm

\par - Let $Z$ be the smallest sub-$\gal$-module of $X^{\du} \otimes Y,$ 
containing $Z_1$ and such that the sub-torus 
$(Z/ Z_1)(1)$ of $(X^{\du} \otimes Y/ Z_1)(1)$ contains
$\pi ({\widetilde b}) $. If we put $Z_2=Z/ Z_1$, we have that 
$Z(1)=Z_1(1) \times Z_2(1). $ 
\vskip 0.3 true cm

\pn With these notation, the Lie algebra $({\Gr}^W_*(W_{-1}\L \G(M)),[\,,\,])$ is 
the Lie algebra $(B+Z(1),[\,,\,])$ and the unipotent radical $W_{-1}(\L \G (M))$ of the Lie algebra of $\G (M)$ is the extension of the abelian variety $B$ by the torus $Z(1)$
defined by the adjoint action of $(B+Z(1),[\,,\,])$ on itself.
\vskip 0.5 true cm

\pn 4.4. {\bf Proposition}
\pn {\it The derived group of the unipotent radical $W_{-1}(\L \G (M))$ of the Lie algebra of $\G (M)$ is the torus $Z_1(1)$.}
\vskip 0.5 true cm

\pn PROOF: This result is a consequence of the motivic version of the structural lemma 
[1] 1.4. In fact, let $g_1=(P,Q,{\vec v}) $ and $
g_2=(R,S,{\vec w} ) $ be two elements of $B+Z(1)$, with $P,R \in B \cap X^{\du} \otimes A(k) $, $ Q,S \in B \cap A^* \otimes Y(k)$, and 
${\vec v},{\vec w} \in Z(1)(k).$ 
We have 

$$g_1 \circ g_2=( P+R,Q+S,{\vec v} +{\vec w} +\Upsilon (P,Q,R,S)
\leqno(4.4.1)$$

\pn where $\Upsilon$ is a $\gal$-equivariant homomorphism from 
$ (X^{\du} \otimes A+A^* \otimes Y)(\ok)\otimes (X^{\du} \otimes A+A^* \otimes Y)(\ok)$
 to $ X^{\du} \otimes Y(\ok)$.
In order to determine $\Upsilon$, we have to understand how $g_1 \circ g_2$ acts on 
the 1-motive $M$. 
The 1-motives $M/W_{-2}M$ and $W_{-1}M$ generate Tannakian subcategories of 
$\langle M/W_{-2}M + W_{-1}M \rangle^\otimes $ and therefore by 1.5 we have
 the surjective morphisms 

$$\eqalign{
pr_1: W_{-1} \L \G( \langle M/W_{-2}M + W_{-1}M \rangle^\otimes ) &\longrightarrow 
W_{-1} \L \G( \langle M/W_{-2}M \rangle^\otimes) \cr
pr_2: W_{-1} \L \G( \langle M/W_{-2}M + W_{-1}M \rangle^\otimes )&\longrightarrow 
W_{-1} \L \G(\langle W_{-1}M \rangle^\otimes ).\cr}\leqno(4.4.2)$$

\pn Since by [2] 3.10 ${\Gr}^W_{-1} \L \G(M)$ is the abelian variety $B$ and  $W_{-2} \L \G(M)$ is the torus $Z(1)$, according the third short exact sequence of (4.2.2) we get that $W_{-1} \L \G( \langle M/W_{-2}M + W_{-1}M \rangle^\otimes )$ is the abelian variety $B$. It follows that explicitly the morphisms (4.4.2) are the projections

$$\eqalign{
pr_1: B &\longrightarrow B \cap X^{\du} \otimes A\cr
pr_2: B &\longrightarrow B \cap A^*\otimes Y.\cr}$$

\pn Let $\pi: W_{-1}(\L \G(M)) \longrightarrow  W_{-1}( \L\G (\langle M/W_{-2}M + W_{-1}M \rangle^\otimes) )$ be the surjective morphism coming from the faithfully flat morphism of the third short exact sequence of (4.2.2).
By definition of $W_{-1}(\L\G (M/W_{-2}M))$ we have that

$$\eqalign{ \big(pr_1(\pi \,g_1)-id \big) \,  W_0/W_{-2}(M) &
\subseteq A,\cr
 \big(pr_1(\pi \,g_1)-id \big) \, A&=0.\cr}$$

\pn Hence modulo the canonical isomorphism 
 ${\underline {\Hom}}(X;A) \cong  X^{\du} \otimes A$ 
 which allows us to identify
$pr_1(\pi \,g_1)-id$ with $P \in B \cap X^{\du} \otimes A(k)$, we obtain that

$$P:{\rm Gr}_{0}^{W}(M) \longrightarrow A.\leqno(4.4.3)$$

\pn In an analogous way, by definition of $W_{-1}(\L\G (W_{-1}M))$ we observe
 that

$$\eqalign{ 
\big(pr_2(\pi \,g_2)-id \big) \,W_{-1}M & \subseteq T,\cr
\big(pr_2(\pi \,g_2)-id \big) \, T& =0.\cr}\leqno(4.4.4)$$

\pn Since the Cartier dual of $W_{-1}(M)$ is the 1-motive $M^*/ W_{-2}M^*,$
 $ pr_2( \pi \,g_2)$ acts on a contravariant way on 
 $M^*/ W_{-2}M^*,$ and therefore we have that

$$\eqalign{ 
\big(pr_2(\pi \,g_2)'-id \big) \, M^*/ W_{-2}M^* & \subseteq A^*,\cr
\big(pr_2(\pi \,g_2)'-id \big) \, A^*& =0,\cr}$$

\pn where the symbol $'$ denote the contravariant action. Consequently,
 modulo the canonical isomorphism $ {\underline {\Hom}}(Y^{\du};A^*) \cong A^*
 \otimes Y$ which allows us to identify
\pn $pr_2(\pi \,g_2)'-id$ with $S \in B \cap A^* \otimes Y$, we have by (4.4.4) that

$$-S:A  \longrightarrow Y(1) .\leqno(4.4.5)$$

\pn Again modulo the canonical isomorphism ${\underline {\Hom}}(X;Y(1)) \cong 
 X^{\du} \otimes Y(1),$ from (4.4.3) and (4.4.5) we get

$$ -[P;S]: {\rm Gr}_{0}^{W}( M)
{\buildrel P \over \longrightarrow}  A  {\buildrel -S \over \longrightarrow} Y(1).  $$

\pn It follows that 
 $\Upsilon (P,Q,R,S)=-[P;S]$ and using (4.4.1) we can conclude that 

$$g_1 \circ g_2 - g_2 \circ g_1=( 0,0,-[P;S]+[R;Q])$$

\pn which is an element of $Z_1(1)$ by definition.
\vskip 0.5 true cm

\pn 4.5. Let $\{e_i\}_i$ and $\{f_j^*\}_j$ be basis of
 $X(\ok)$ and $Y^{\du}(\ok)$ respectively. 
Choose a point $P$ of $B \cap X^{\du} \otimes A (k)$ and a point
$Q$ of $B \cap A^* \otimes Y (k)$ such that the abelian sub-variety they generate in 
$X^{\du} \otimes A +A^* \otimes Y$, is isogeneous to  $B$. 
Denote by ${\overline v}:X(\ok) \longrightarrow A(\ok)$ et 
${\overline v}^*:Y^{\du}(\ok) \longrightarrow 
A^*(\ok)$ the $\gal$-equivariant homomorphisms defined by 

$$\eqalign{{\overline v}(e_i)&=\alpha(P,e_i),\cr 
{\overline v}^*(f_j^*)&=\beta(Q,f_j^*).\cr}$$

\pn Moreover choose a point
${\vec q}=(q_1,\dots,q_{{\rm rg}\, Z_2})$ 
of $Z_2(1)(k)$ such that the points
$q_1, \dots,q_{{\rm rg}\, Z_2} $
 are multiplicative independent. 

\pn Let $ \Gamma: Z(1)(\ok) \otimes X \otimes Y^{\du}(\ok)
\longrightarrow  \ZZ(1)(\ok)$ be the $\gal$-equivariant homomorphism obtained from the maps 
$ \gamma: (X^{\du} \otimes Y(1)) \otimes X 
\longrightarrow  Y(1)$  and
$ev_Y: Y \otimes Y^{\du} \longrightarrow \ZZ(0)$, and denote by  
${\overline \psi}: X \otimes Y^{\du}(\ok) \longrightarrow {\ZZ}(1)(\ok)$ 
the $\gal$-equivariant homomorphism defined by 

$${\overline \psi}(e_i,f_j^*)= \Gamma([P,Q],{\vec q},e_i,f_j^* ).\leqno(4.5.1)$$ 
\vskip 0.5 true cm

\pn 4.6. {\bf Lemma}
\pn {\it With the above notation, the Tannakian category $\langle M \rangle^\otimes$ is equivalent to the Tannakian category generated by the 1-motives 
$(e_i\ZZ,f_j^*\ZZ,A,A^*,{\overline v}_{i},{\overline v}^*_{j}, {\overline \psi}_{i,j})$, where ${\overline v}_{i},{\overline v}^*_{j}$ and $ {\overline \psi}_{i,j}$ are the $\gal$-equivariant homomorphisms obtained restricting 
respectively ${\overline v},{\overline v}^*$ and $ {\overline \psi}$  to 
$e_i\ZZ$ and $f_j^*\ZZ$: }

$$\eqalign{
{\overline v}_{i}: e_i\ZZ \longrightarrow A,& \qquad  e_i \longmapsto \alpha(P,e_i),\cr
{\overline v}^*_{j}:f_j^* \ZZ \longrightarrow A^*, & \qquad  
f_j^*  \longmapsto \alpha^*(Q,f_j^*),\cr
{\overline \psi}_{i,j}:e_i\ZZ \times f_j^*\ZZ  \longrightarrow {\cal P}_{\vert e_i\ZZ \times f_j^*\ZZ },& \qquad (e_i,f_j^*) \longmapsto \Gamma([P,Q],{\vec q},e_i,f_j^* ).\cr}$$

\pn PROOF: According to the proof of theorem 3.8 [2], the homomorphisms ${\overline v},
{\overline v}^*$ and $ {\overline \psi}$  define a 1-motive 
$(X,Y^{\du},A,A^*,{\overline v},{\overline v}^*, {\overline \psi})$ which
 generates the same Tannakian category as $M$. In order to conclude we apply theorem
 [1] 1.7  to the 1-motives 
\pn $(X,Y^{\du},A,A^*,{\overline v},{\overline v}^*, {\overline \psi})$
 and $(e_i\ZZ,f_j^*\ZZ,A,A^*,{\overline v}_{i},{\overline v}^*_{j}, {\overline \psi}_{i,j})$. 
\vskip 0.5 true cm

\pn 4.7. Consider the 1-motives 

$$\eqalign{
M^{tab}=&\oplus_{i,j}
 (e_i\ZZ,f_j^*\ZZ,0,0,0,0,{\overline \psi}_{i,j}^{ab})\cr
M^{a}=&\oplus_{i,j}
 (e_i\ZZ,f_j^*\ZZ,A,A^*,{\overline v}_{i},{\overline v}^{*}_{j},0)\cr
M^{nab}=&\oplus_{i,j}
 (e_i\ZZ,f_j^*\ZZ,A,A^*,{\overline v}_{i},{\overline v}^{*}_{j},
{\overline \psi}_{i,j}^{nab})\cr
M^{ab}=&\oplus_{i,j}
 (e_i\ZZ,f_j^*\ZZ,A,A^*,{\overline v}_{i},{\overline v}^{*}_{j},
{\overline \psi}_{i,j}^{ab})\cr}$$

\pn where

$$\eqalign{ 
{\overline \psi}_{i,j}^{nab}(e_i,f_j^*)&=\Gamma([P,Q],{\vec 1},e_i,f_j^*)\cr
{\overline \psi}_{i,j}^{ab}(e_i,f_j^*)&=\Gamma({\vec 1},{\vec q},e_i,f_j^*)\cr
}\leqno(4.7.1)$$

\pn According to lemma 4.6, the 1-motives $M^{tab},M^{a},M^{ab}$ and $M^{nab}$ belong 
to the Tannakian category $\langle M \rangle^\otimes $ generated by $M$.  
\vskip 0.5 true cm

\pn 4.8. {\bf Lemma}
\pn {\it The Tannakian category generated by $M$ is equivalent to the Tannakian 
category generated by the 1-motive $ M^{tab} \oplus M^{nab}$. Moreover the 1-motives
$M^{ab}$ and $M^{a} \oplus M^{tab}$ generate the same Tannakian category.}
\vskip 0.5 true cm 

\pn PROOF: Since through the projection
 $p:X^{\du} \otimes Y(1) \longrightarrow (X^{\du} \otimes Y/ Z_1)(1)$  
the $\Sigma-(X^{\du} \otimes Y)(1)$-torsor $i^* {\cal B}$ becomes a trivial 
torsor, i.e. $p_*i^* {\cal B}= B \times (X^{\du} \otimes Y/ Z_1)(1)$,
confronting (4.5.1) and (4.7.1), we observe that
to have the trivialization
$ {\overline \psi}_{i,j}^{ab}$ and ${\overline \psi}_{i,j}^{nab}$ 
is the same thing as to have the trivialization 
$ {\overline \psi}_{i,j}$. Hence by lemma 4.6, 
we have that the 1-motives  $M$ and $ M^{tab} \oplus M^{nab}$ generate 
the same Tannakian category.
\pn Always because of the fact that the  $\Sigma-(X^{\du} \otimes Y/Z_1)(1)$-torsor
$p_*i^* {\cal B}$ is trivial, we observe that the trivialization 
${\overline \psi}_{i,j}^{ab}$ is independent of the abelian part 
of the 1-motive M, i.e. it is independent of ${\overline v}_{i},{\overline v}^{*}_{j}$ Therefore, we can conclude that the Tannakian category generated by 1-motive $M^{ab}$ is equivalent to the Tannakain category generated by the 1-motive $M^{a} \oplus M^{tab}.$  
\vskip 0.5 true cm

\pn 4.9. {\bf Theorem}
\pn {\it The Tannakian category generated by $M^{ab}$ 
is the biggest Tannakian subcategory of $\langle M \rangle^\otimes $ whose motivic Galois group is commutative.
We have the following diagram of morphisms of 
affine group sub-$\langle M \rangle^\otimes $-schemes}

$$\matrix{
{\scriptstyle 0} &\rightarrow &{\scriptstyle Z_{2}(1)} &\rightarrow &{\scriptstyle \L\G(M)} & {\buildrel i^{nab} \over  \rightarrow } &{\scriptstyle I^{nab}\L\G(M^{nab})} &\rightarrow &{\scriptstyle 0} \cr
&&\downarrow &&||&&\downarrow& \cr
{\scriptstyle 0} &\rightarrow &{\scriptstyle Z(1)} &\rightarrow &{\scriptstyle \L\G(M)} & {\buildrel i^{a}\over  \rightarrow } &{\scriptstyle I^{a} \L\G(M^{a})} &\rightarrow &{\scriptstyle 0} \cr
&&\uparrow &&||&&\uparrow& \cr
{\scriptstyle 0} &\rightarrow &{\scriptstyle Z_1(1)} &\rightarrow &{\scriptstyle \L\G(M)} & {\buildrel i^{ab}\over  \rightarrow } &{\scriptstyle I^{ab}\L\G( M^{ab})} &\rightarrow &{\scriptstyle 0} \cr
&&\downarrow &&||&&\downarrow& \cr
{\scriptstyle 0} &\rightarrow &{\scriptstyle B+Z_1(1)} &\rightarrow &{\scriptstyle \L\G(M)} & {\buildrel i^{tab}\over  \rightarrow } &{\scriptstyle I^{tab} \L\G(M^{tab})} &\rightarrow &{\scriptstyle 0} \cr}$$

\pn {\it where all horizontal short sequences are exact and where the vertical arrows on the left are inclusions and those on the right are surjections.}
\vskip 0.5 true cm

\pn PROOF: By [2] 3.10 we know that ${\Gr}^W_* (W_{-1} \L \G(M))$ is the Lie algebra $(B+Z(1),[ \, ,\,])$. Now, from the definition of the 1-motives
 $ M^{nab}, M^{ab}$
and $M^{tab}$, and from (4.7.1) we observe that the torus $Z_2(1)$ acts trivially on $ M^{nab}$, that the torus $Z_1(1)$ acts trivially on $ M^{ab}$ and that
the split 1-motive $B+Z_1(1)$ acts trivially on $M^{tab}.$ In other words we have that 

$$\eqalign{
&\L H_{\langle M \rangle^\otimes}(\langle M^{nab} \rangle^\otimes)=Z_2(1)\cr
&\L H_{\langle M \rangle^\otimes}(\langle M^{ab} \rangle^\otimes)=Z_1(1)\cr
& \L H_{\langle M \rangle^\otimes}(\langle M^{tab} \rangle^\otimes)=B+Z_1(1)\cr}$$

\pn and therefore we get the first, the second and the fourth short exact sequence.

\pn Recall that by [2] 3.10 the Lie algebra $W_{-2}\L \G (M)$ is the torus $Z(1)$.
Moreover by construction, the 1-motive without toric part $M^{a}$ generates the same
 Tannakian category as the 1-motive $W_0/W_{-2}M+W_{-1}M.$  Hence thanks to (4.2.1) the obtain the second horizontal short exact sequence. ( It is possible to prove the exactness of the second short sequence computing the Lie algebra $\L 
H_{\langle M \rangle^\otimes}(\langle M^{a} \rangle^\otimes).)$

\par According to lemma 4.8, the 1-motives $M^{a}$ and $M^{tab}$ generate Tannakian subcategories of the Tannakian category $\langle M^{ab} \rangle^\otimes.$
By construction the 1-motive $M^{a}$ generates a Tannakian subcategory of 
the Tannakian category $M^{nab}$. Hence in order to prove that the left vertical arrows are inclusions and that 
the right vertical arrows are surjections, it is enough to apply the lemma 1.8.

\par The third exact sequence of the above diagram implies that the motivic Galois group of $M^{ab}$ is isomorph to the quotient 

$$ \L \G (M) / Z_1(1)$$

\pn But according proposition 4.4, $Z_1(1)$ is the derived group of $\L \G (M)$
 and hence we can conclude that $M^{ab}$ generate the biggest Tannakian
 subcategory of 
$\langle M \rangle^\otimes$ whose motivic Galois group is commutative.  
\vskip 0.5 true cm

\pn 4.10. REMARK: Among the non degenerate 1-motives, the 1-motive $M^{nab}$ is the one which generates the biggest Tannakian subcategory of $\langle M \rangle^\otimes$, whose motivic Galois group is non commutative. (A 1-motive is said to be non degenerate if the dimension of the Lie algebra $W_{-1} \L \G(M)$ is  maximal (cf. [1] 2.3)). 
\vskip 0.5 true cm

\pn {\bf Bibliography}
\vskip 0.5 true cm

\par\noindent
1. C. Bertolin, {\it The Mumford-Tate group of 1-motives,} Ann. Inst. Fourier, 52, 4 (2002).

\par\noindent
2. C. Bertolin, {\it Le radical unipotent du groupe de Galois motivique d'un 1-motif,} Math. Ann. 327 (2003).

\par\noindent
3. J.-L. Brylinski, {\it 1-motifs et formes automorphes
 (th\'eorie arithm\'etique des domaines 
des Siegel)}, Pub. Math. Univ. Paris VII 15 (1983).

\par\noindent
4. P. Deligne, {\it Th\' eorie de Hodge III}, Pub. Math. de
 l'I.H.E.S. 44 (1975).

\par\noindent
5. P. Deligne and J.S. Milne, {\it Tannakian Categories}, LN 900 (1982).

\par\noindent
6.  P. Deligne, {\it Le groupe fondamental de la droite 
projective moins trois points}, Galois group over $\QQ$, Math. Sci. Res. Inst. Pub. 16 (1989).

\par\noindent
7. P. Deligne, {\it Cat\'egories tannakiennes}, The Grothendieck Festschrift II,
Birk\-h\"au\-ser 87 (1990).

\par\noindent
8. U. Jannsen, {\it Mixed motives and algebraic K-theory}, LN 1400 (1990).

\par\noindent
9. U. Jannsen, {\it Mixed motives, motivic cohomology, and Ext-groups}, 
Proceedings of the ICM (Z\"urich, 1994), Birkh\"auser (1995).

\par\noindent
10. J.S. Milne, {\it Motives over finite fields}, Motives, Proc. of Symp. in Pure Math. 55 (1994).

\par\noindent
11. N. Saavedra Rivano, {\it Cat\' egories tannakiennes}, LN 265 (1972).

\vskip 2 true cm

\pn \hskip 1 true cm Bertolin Cristiana
\pn \hskip 1 true cm NWF I - Mathematik
\pn \hskip 1 true cm Universit\"at Regensburg
\pn \hskip 1 true cm D-93040 Regensburg
\pn \hskip 1 true cm Germany
\vskip 0.2 true cm
\pn \hskip 1 true cm Email: cristiana.bertolin@mathematik.uni-regensburg.de

\end